\documentclass[12pt,oneside]{article}
\usepackage{amsmath,amssymb,amsfonts,amsthm}
\usepackage{color}

\newcommand{\T}{{\cal T}}

\newcommand{\Real}{\mathbb R}

\newcommand{\To}{\longrightarrow}

\newcommand{\prof}{\noindent \textit{\textbf{Proof.\:\:}}}

\newcommand{\p}{\pi^{-1}(TM)}

\def\o#1{\overline{#1}}

\def\pa{\partial}
\def\paa{\dot{\partial}}

\def\Section#1{\vspace{30truept}\addtocounter{section}{1}\setcounter{thm}{0}\setcounter{equation}{0}
{\noindent\Large\bf\arabic{section}.~~#1}\par \vspace{12pt}}

\newtheorem{thm}{Theorem}[section]
\newtheorem{cor}[thm]{Corollary}
\newtheorem{lem}[thm]{Lemma}
\newtheorem{prop}[thm]{Proposition}
\newtheorem{defn}[thm]{Definition}

\newtheorem{rem}[thm]{Remark}

\numberwithin{equation}{section}

\title{\bf A GLOBAL THEORY OF CONFORMAL FINSLER GEOMETRY}

\author{\textbf{Nabil L. Youssef$\,\dagger$ , S. H. Abed$\dagger$ }\\
\bf{and A. Soleiman$\ddagger$}}
\date{}

\begin{document}               
\bibliographystyle{plain}
\maketitle                     
\vspace{-1.15cm}
\begin{center}
{$\dagger$Department of Mathematics, Faculty of Science,\\ Cairo
University, Giza, Egypt.}
\end{center}
\vspace{-0.8cm}
\begin{center}
nyoussef@frcu.eun.eg,\ sabed@frcu.eun.eg
\end{center}
\vspace{-0.7cm}
\begin{center}
and
\end{center}
\vspace{-0.7cm}
\begin{center}
{$\ddagger$Department of Mathematics, Faculty of Science,\\ Benha
University, Benha, Egypt.}
\end{center}
\vspace{-0.8cm}
\begin{center}
soleiman@mailer.eun.eg
\end{center}
\smallskip
\begin{center}
{\textbf{Dedicated to the memory of Prof. Dr. A. TAMIM}}
\end{center}

\vspace{1cm} \maketitle
\smallskip

\noindent{\bf Abstract.} The aim of the present paper is to
establish \emph{a global investigation} of conformal changes in
Finsler geometry. Under this change, we obtain the relationships
between some geometric objects associated to $(M,L)$ and the
corresponding objects associated to $(M,\tilde{L})$,
$\widetilde{L}=e^{\sigma(x)} L$ being the Finsler conformal
transformation. We have found explicit global expressions relating
the two associated Cartan connections $\nabla$ and
$\widetilde{\nabla}$, the two associated Berwald connections $D$ and
$\widetilde{D}$ and the two associated Barthel connections $\Gamma$
and $\widetilde{\Gamma}$. The relationships between the
corresponding curvature tensors have been also found. The relations
thus obtained lead in turn to several interesting results.
\par
Among the results obtained, is a characterization of conformal
changes, a characterization of homotheties, some conformal
invariants and conformal $\sigma$-invariants. In addition, several
useful identities have been found.
\par
Although our treatment is entirely global, the local expressions of
the obtained results, when calculated, coincide with the existing
classical local results. \footnote{This paper was presented in \lq
\lq \,The 9 th. International Conference of Tensor Society on
Differential Geometry, informatics and their Applications \rq\rq\
held at Sapporo, Japan, September 4-8, 2006.}



\newpage
\vspace{30truept}\centerline{\Large\bf{Introduction}}\vspace{12pt}
\par Studying Finsler geometry one encounters substantial
difficulties trying to seek analogues of classical global, or
sometimes even local, results of Riemannian geometry. These
difficulties arise mainly from the fact that in Finsler geometry all
geometric objects  depend not only on positional coordinates, as in
Riemannian geometry, but also on directional arguments.

 In Riemannian geometry there is a canonical linear connection on the
manifold $M$, whereas in Finsler geometry there is a corresponding
canonical linear connection due to E. Cartan. However, this is not a
connection on $M$ but in $T(\T M)$, the tangent bundle of $\,\T M$
(\emph{the Klein-Grifone approach}), or in $\,\pi^{-1}(TM) $, the
pullback of the tangent bundle $TM$ by $\,\pi: \T M\longrightarrow
M$ (\emph{the pullback approach}).
\par
The infinitesimal transformations in Finsler geometry (such as
conformal~\cite{conf. 3}, projective~\cite{del castillo},
semi-projective~\cite{Nabil.1}, $\beta-$changes~\cite{Shimada},
conformal $\beta$-changes~\cite{Abed}, ... etc.) play an important
role not only in differential geometry but also in application to
other branches of science, especially in the process of
geometrization of physical theories. Conformal changes have been
initiated  by M. S. Kneblman~\cite{conf. 1} and have been
investigated  by many authors \cite{conf. 3},\,\cite{conf.
2},\,\cite{conf. 4},\,\cite{Miron 2},...etc. Almost all known result
concerning these changes are local ones. The global results are very
few in the literature. In this paper we present \emph{a global
theory} of conformal changes in Finsler geometry

\par
The most well-known and widely used approaches to GLOBAL Finsler
geometry are the Klein-Grifone (KG-) approach (cf. \cite{Grifone
1},~\cite{Grifone 2},~\cite{Klien1} and~\cite{Nabil.1}) and the
pull-back \linebreak(PB-) approach (cf.~\cite{Akbar 1},\,\cite{Akbar
2},\,\cite{Dazord},\,\cite{Mats1} and~\,\cite{Ali1}). The universe
of the first approach is the vector bundle
$\,\pi_{TM}:TTM\longrightarrow TM$, whereas the universe of the
second is the vector bundle $P:\pi^{-1}(TM)\longrightarrow TM$. Each
of the two approaches has its own geometry which differs
significantly from the geometry of the other (in spite of the
existence of some links between them). Each also has its advantages
and disadvantages. For example, the KG-formalism is an elegant one
and easily manipulated, but the geometric aspects of many of its
objects are not clarified. On the other hand, the PB-formalism is
similar to and guided by Riemannian geometry, but is somewhat
difficult to manipulate.
\par
Most of the geometers, when treating Finsler geometry from a GLOBAL
standpoint, follow \textbf{uniquely} one of the above mentioned
approaches. We proceed here differently. We establish a global
theory of conformal Finsler geometry within the PB-approach, making
simultaneous use of the KG-approach. This has been done via certain
links we have found between both approaches. This \lq\lq\,double
approach" enables us to overcome several difficulties and to
continue our development. For instance, without the insertion of the
KG-approach, we were unable to find the relation between the Cartan
connection and its conformal transform, a relation which is
fundamental for the present work. This shows that these two
approaches are not alternatives but rather complementary.
\par
Let $(M, L)$ and $(M, \widetilde{L})$ be two conformal Finsler
manifolds, where the conformal transformation is given by
$L\longrightarrow\widetilde{L}=e^{\sigma(x)} L$. We have found
explicit global expressions relating the two associated Cartan
connections $\nabla$ and $\widetilde{\nabla}$, the two associated
Berwald connections $D$ and $\widetilde{D}$ and the two associated
Barthel connections $\Gamma$ and $\widetilde{\Gamma}$. The
relationships between the corresponding curvature tensors have also
been found. The relations thus obtained lead in turn to several
interesting results.
\par
Among the results obtained, is a characterization of conformal
changes, a characterization of homotheties, some conformal
invariants and conformal $\sigma$-invariants. In addition, several
useful identities have been found.
\par
The last section of the paper presents an application of some of the
results obtained. It deals with geodesics and Jacobi fields in the
context of global Finsler geometry.
\par
It should finally be noted that, although our treatment is entirely
global, the local expressions of the obtained results, when
calculated, coincide with the classical local results of Hashiguchi,
Matsumoto and others. For the sake of completeness an appendix,
concerning the local expressions of the most important geometric
objects treated, is included.


\Section{Fundamentals of the pull-back formalism}

In this section we give a brief account of the basic concepts
 of the pullback formalism necessary for this work. For more details
  refer to~\cite{Akbar 1},\,\cite{Akbar
2},\,\cite{Dazord},\,\cite{Mats1} and~\,\cite{Ali1}.
 We make the general
assumption that all geometric objects we consider are of class
$C^{\infty}$. The
following notation will be used throughout this paper:\\
 $M$: a real differentiable manifold of finite dimension $n$ and of
class $C^{\infty}$,\\
 $\mathfrak{F}(M)$: the $\Real$-algebra of differentiable functions
on $M$,\\
 $\mathfrak{X}(M)$: the $\mathfrak{F}(M)$-module of vector fields
on $M$,\\
$\pi_{M}:TM\longrightarrow M$: the tangent bundle of $M$,\\
$\pi: \T M\longrightarrow M$: the subbundle of nonzero vectors
tangent to $M$,\\
$V(TM)$: the vertical subbundle of the bundle $TTM$,\\
 $P:\pi^{-1}(TM)\longrightarrow \T M$ : the pullback of the
tangent bundle $TM$ by $\pi$,\\
 $\mathfrak{X}(\pi (M))$: the $\mathfrak{F}(M)$-module of
differentiable sections of  $\pi^{-1}(T M)$,\\
$ i_{X}$ : interior product with respect to  $X
\in\mathfrak{X}(M)$,\\
$df$ : the exterior derivative  of $f$,\\
$ d_{L}:=[i_{L},d]$, $i_{L}$ being the interior derivative with
respect to the vector form $L$.
\par Elements  of  $\mathfrak{X}(\pi (M))$ will be called
$\pi$-vector fields and will be denoted by barred letters
$\overline{X} $. Tensor fields on $\pi^{-1}(TM)$ will be called
$\pi$-tensor fields. The fundamental $\pi$-vector field is the
$\pi$-vector field $\overline{\eta}$ defined by
$\overline{\eta}(u)=(u,u)$ for all $u\in \T M$. The lift to
$\pi^{-1}(TM)$ of a vector field $X$ on $M$ is the $\pi$-vector
field $\overline{X}$ defined by $\overline{X}(u)=(u,X(\pi (u)))$.
The lift to $\pi^{-1}(TM)$ of a $1$-form $\omega$ on $M$ is the
$\pi$-form $\overline{\omega}$ defined by
$\overline{\omega}(u)=(u,\omega(\pi (u)))$.

 The tangent bundle $T(\T M)$
is related to the pullback bundle $\pi^{-1}(TM)$ by the short exact
sequence \vspace{-0.4cm}
$$0\longrightarrow
 \pi^{-1}(TM)\stackrel{\gamma}\longrightarrow T(\T M)\stackrel{\rho}\longrightarrow
\pi^{-1}(TM)\longrightarrow 0 ,\vspace{-0.2cm}$$
 where the bundle morphisms $\rho$ and $\gamma$ are defined respectively by
$\rho = (\pi_{\T M},d\pi)$ and $\gamma (u,v)=j_{u}(v)$, where
$j_{u}$  is the natural isomorphism $j_{u}:T_{\pi_{M}(v)}M
\longrightarrow T_{u}(T_{\pi_{M}(v)}M)$. The vector $1$-form $J$ on
$TM$ defined by $J=\gamma o \rho$ is called the natural almost
tangent structure of $T M$. The vertical vector field $\mathcal{C}$
on $TM$ defined by $\mathcal{C}:=\gamma o \overline{\eta} $ is
called the fundamental or the canonical (Liouville) vector field.

 Let $\nabla$ be  a linear
connection (or simply a connection) in the pullback bundle
$\pi^{-1}(TM)$.
 We associate to
$\nabla$ the map\vspace{-0.2cm}
$$K:T \T M\longrightarrow \pi^{-1}(TM):X\longmapsto \nabla_X \overline{\eta}
,\vspace{-0.2cm}$$ called the connection (or the deflection) map of
$\nabla$. A tangent vector $X\in T_u (\T M)$ is said to be
horizontal if $K(X)=0$ . The vector space $H_u (\T M)= \{ X \in T_u
(\T M) : K(X)=0 \}$ of the horizontal vectors
 at $u \in  \T M$ is called the horizontal space to $M$ at $u$  .
   The connection $\nabla$ is said to be regular if
 $$T_u (\T M)=V_u (\T M)\oplus H_u (\T M) \qquad \forall u\in \T M .$$
  \par If $M$ is endowed with a regular connection, then the vector bundle
   maps
\begin{eqnarray*}
 \gamma &:& \pi^{-1}(T M)  \To V(\T M), \\
   \rho |_{H(\T M)}&:&H(\T M) \To \pi^{-1}(TM), \\
   K |_{V(\T M)}&:&V(\T M) \To \pi^{-1}(T M)
\end{eqnarray*}
 are vector bundle isomorphisms.
   Let us denote
 $\beta=(\rho |_{H(\T M)})^{-1}$,
then \vspace{-0.2cm}
   \begin{align}\label{fh1}
    \rho  o  \beta = id_{\pi^{-1} (\T M)}, \quad  \quad
       \beta o \rho =\left\{
                                \begin{array}{ll}
                                          id_{H(\T M)} & \hbox{on   H(\T M)} \\
                                         0 & \hbox{on    V(\T M)}
                                       \end{array}
                                     \right.\vspace{-0.2cm}
\end{align}

The classical  torsion tensor $\textbf{T}$  of the connection
$\nabla$ is defined by
$$\textbf{T}(X,Y)=\nabla_X \rho Y-\nabla_Y\rho X -\rho [X,Y] \quad
\forall\,X,Y\in \mathfrak{X} (\T M)$$ The horizontal and mixed
torsion tensors, denoted respectively by $Q $ and $ T $, are defined
by \vspace{-0.2cm}
$$Q (\overline{X},\overline{Y})=\textbf{T}(\beta \overline{X}\beta \overline{Y}),
\, \,\, T(\overline{X},\overline{Y})=\textbf{T}(\gamma
\overline{X},\beta \overline{Y}) \quad \forall \,
\overline{X},\overline{Y}\in\mathfrak{X} (\pi (M)).\vspace{-0.2cm}$$
\par The classical curvature tensor  $\textbf{K}$ of the connection
$\nabla$ is defined by
 $$ \textbf{K}(X,Y)\rho Z=-\nabla_X \nabla_Y \rho Z+\nabla_Y \nabla_X \rho Z+\nabla_{[X,Y]}\rho Z
  \quad \forall\, X,Y, Z \in \mathfrak{X} (\T M).$$
The horizontal, mixed and vertical curvature tensors, denoted
respectively by $R$, $P$ and $S$, are defined by
$$R(\overline{X},\overline{Y})\o Z=\textbf{K}(\beta
\overline{X}\beta \overline{Y})\o Z,\quad
P(\overline{X},\overline{Y})\o Z=\textbf{K}(\beta
\overline{X},\gamma \overline{Y})\o Z,\quad
S(\overline{X},\overline{Y})\o Z=\textbf{K}(\gamma
\overline{X},\gamma \overline{Y})\o Z.$$

\begin{defn}A Finsler manifold of dimension $n$ is a pair $(M,L)$, where $M$ is a
 differentiable manifold of dimension $n$ and $L: TM \To \Real
 $ is a map {\em ({\it called Lagrangian})} such that\,{\em:}
 \begin{description}
   \item[(a)] $L(X)\geqq 0,$\, for all $X \in \T M $,  $L(0)=0 $,
   \item[(b)] $L $ is  $C^{\infty}$ on  $\T M$ ,  $C^{1}$ on $TM$,
   \item[(c)]$L$ is homogenous of degree 1 in the directional argument $y$\,{\em:} $\mathcal{C} \cdot L=L$ ,
   \item[(d)] The $\pi$-form $g$ of order two with components
   $g_{ij}:= \frac{1}{2}\dot{\partial_{j}} \dot{\partial_{i}} L^2$ is positive
   definite. \vspace{-0.2cm}
 \end{description}
 \end{defn}
 \noindent{\it Hence, the $\pi-$form $g$ defines a positive definite metric in
   $\pi^{-1}(TM)$}\\
 \noindent ({ Here, $\dot{\partial}_{i}$ denotes partial differentiation with respect to the directional
   argument $y^i$}).

\begin{thm} \label{th.1}{\em{\cite{Ali2}}} Let $(M,L)$ be a Finsler manifold. There exists a
unique regular connection $\nabla$ in $\pi^{-1}(TM)$ such that
\begin{description}
  \item[(a)]  $\nabla$ is  metric\,{\em:} $\nabla g=0$,

  \item[(b)]   The horizontal torsion of $\nabla$ vanishes\,{\em:} $Q=0
  $,
  \item[(c)]  The mixed torsion $T$ of $\nabla$ satisfies \,
  $g(T(\overline{X},\overline{Y}), \overline{Z})=g(T(\overline{X},\overline{Z}),\overline{Y})$.
\end{description}
\end{thm}
Such a connection is called the Cartan connection associated to the
Finsler manifold $(M,L)$.

 For the Cartan connection $\nabla$, we have
\vspace{-0.2cm}
   \begin{align}\label{fh2}
    Ko\gamma = id_{\pi^{-1} (\T M)}, \quad \quad
      \gamma o K =\left\{
                                       \begin{array}{ll}
                                          id_{V(\T M)} & \hbox{on   V(\T M)} \\
                                         0 & \hbox{on    H(\T M)}
                                       \end{array}
                                     \right.
\end{align}
   Then, from (\ref{fh1}) and  (\ref{fh2}), we get\vspace{-0.2cm}
 \begin{equation}\label{hv}
   \beta o \rho + \gamma o K = id_{T(\T M)}\vspace{-0.2cm}
\end{equation}
 Hence, if we set
   $h {\!\!}:= \beta o \rho $, $v {\!\!}:=\gamma o K$, then every vector field $X{\!\!}\in{\!\!}\mathfrak{X}(\T M)$ can be
   represented uniquely in the form\vspace{-0.2cm}
           \begin{equation}\label{uy}
                                         X = hX +vX =\beta \rho X + \gamma K X \vspace{-0.2cm}
                                  \end{equation}
 The maps $h$ and $v$ are the horizontal and  vertical projectors associated
 to the Cartan connection $\nabla$: $h^{2}= h$, \,\,\,$v^{2}= v$,\,\,\,
  $h+v= id_{\mathfrak{X}(TM)}$,\,\, \, $voh=hov=0 .$

One can show that the torsion  of the Cartan connection has the
property that $T(\overline{X},\overline{\eta})=0$ for all
$\overline{X} \in \mathfrak{X} (\pi (M))$, and associated to $T$ we
have the

\begin{defn}{\em{\cite{Ali2}}} Let  $\nabla$ be the Cartan  connection associated to $(M,L)$.
The torsion tensor field  $T$ of the connection $\nabla$ induces a
$\pi$-tensor field of type $(0,3)$,  denoted again  $T$,  defined
by{\,\rm :}
$$T(\overline{X},\overline{Y},\overline{Z})=g(T(\overline{X},\overline{Y}),\overline{Z}),
\quad \text {for all} \quad
\overline{X},\overline{Y},\overline{Z}\in \mathfrak{X} (\pi (M))$$
and a $\pi$-form $C$ defined by: $$C(\overline{X}):= \text {Trace of
the map} \quad \overline{Y} \longmapsto
T(\overline{X},\overline{Y}),
 \quad \text {for all} \quad \overline{X}\in \mathfrak{X} (\pi (M)).$$
\end{defn}
\begin{defn} {\em{\cite{Ali2}}} With respect to the Cartan  connection $\nabla$ associated to
$(M,L)$, we have
 \begin{description}
    \item[--]The horizontal and vertical Ricci tensors $Ric^h$ and
    $Ric^v$  are defined respectively by{\,\rm:}
 $$Ric^h(\overline{X},\overline{Y}):= \text {Trace of the map} \quad \overline{Z} \longmapsto
R(\overline{X},\overline{Z})\overline{Y}, \quad \text {for all}
\quad \overline{X},\overline{Y}\in \mathfrak{X} (\pi (M)),$$
 $$Ric^v(\overline{X},\overline{Y}):= \text {Trace of the map} \quad  \overline{Z} \longmapsto
S(\overline{X},\overline{Z})\overline{Y}, \quad \text {for all}
\quad \overline{X},\overline{Y}\in \mathfrak{X} (\pi (M)).$$

   \item[--]The horizontal and vertical Ricci  maps $Ric_0^h$ and $Ric_0^v$ are
defined respectively by{\,\rm:}
$$Ric^h(\overline{X},\overline{Y})=g(Ric_0^h(\overline{X}),\overline{Y}),
\quad \text {for all} \quad \overline{X},\overline{Y}\in
\mathfrak{X} (\pi (M)),$$
$$Ric^v(\overline{X},\overline{Y})=g(Ric_0^v(\overline{X}),\overline{Y}),
\quad \text {for all} \quad \overline{X},\overline{Y}\in
\mathfrak{X} (\pi (M)).$$

    \item[--] The horizontal and vertical scalar
curvatures $Sc^h$ , $Sc^v$ are defined respectively
by{\,\rm:}\vspace{-0.2cm}
\begin{equation*}
     Sc^h:=
\text {Tr}( Ric_0^h), \qquad \qquad  Sc^v:= \text {Tr}  (Ric_0^v
),\vspace{-0.3cm}
\end{equation*}
 \end{description}
where $R$ and $S$ are respectively the horizontal and vertical
curvature tensors of $\nabla$.
\end{defn}


\begin{defn} {\em{\cite{Grifone 1}}} A spray  on $M$ is a vector field $X$ on $TM$,
 $C^{\infty}$ on $\T M$, $C^{1}$ on $TM$, such that
\begin{description}
  \item[(a)] $\rho o X = \overline{\eta}$,

  \item[(b)] $X$ is homogenous of degree  $2$ in $y$\,{\em:}   $[\mathcal{C},X]= X $.
\end{description}
\end{defn}

Let  $E:=\frac{1}{2}L^2$ be  the energy of the Lagrangian $L$. One
can show that
\begin{equation}\label{eq.t4}
d_{J}E(X)= g(\rho X, \overline{\eta}), \quad \text{for all} \quad
X\in \mathfrak{X}(\T M),
\end{equation}
and consequently, $g(\overline{\eta}, \overline{\eta})=L^2$. One can
also show that the exterior $2-$form
 $\Omega:=dd_{J}E$ on $TM$ is
nondegenerate. The form $\Omega$ is called the fundamental
form~\cite{Grifone 1}.

\begin{prop}{\em{\cite{Klien1}}} Let $(M,L)$ be a Finsler manifold. The vector field $G\in\mathfrak{X}(TM)$
 determined by $i_{G}\Omega =-dE$ is a spray,
 called the canonical spray associated to the energy $E$.
 \end{prop}
 One can show, in this case, that $G=\beta o \overline{\eta}$, and
 $G$ is thus horizontal with respect to the Cartan connection $\nabla$.

\begin{thm} {\em{\cite{Ali3}}} \label{th.1a} Let $(M,L)$ be a Finsler manifold. There exists  a
unique regular connection $D$ in $\pi^{-1}(TM)$ such that
\begin{description}
  \item[(a)]  $D$ is torsion free,
  \item[(b)]  The canonical spray $G= \beta o \overline{\eta}$ is
  horizontal with respect to $D,$
  \item[(c)] The mixed curvature $P$ of $D$ satisfies\,{\em:} $P\otimes \overline{\eta}=0$.
\end{description}
\end{thm}
Such a connection is called the Berwald connection associated to the
Finsler manifold $(M,L)$.
 \vspace{5pt}
 \par We terminate this section
by some concepts and results concerning the Klein-Grifone approach.
For more details refer to~\cite{Grifone 1},~\cite{Grifone
2},~\cite{Klien1} and~\cite{Nabil.1}. \vspace{-0.2cm}
\begin{defn} {\em{\cite{Grifone 1}}} A nonlinear connection on $M$ is a vector $1$-form $\Gamma$
on $TM$, $C^{\infty}$ on $\T M$, $C^{o}$ on $TM$, such
that\vspace{-0.3cm}$$J \Gamma=J, \quad\quad \Gamma J=-J .$$
\end{defn}

We have to note that $\Gamma$ defines on $TM$ an almost product
structure ($\Gamma^{2}=I$, where $I$ is the identity on $TM$).
 The  horizontal and vertical projectors $h$  and $v$ associated to $\Gamma$ are defined by
   $h:=\frac{1}{2} (I+\Gamma),\, v:=\frac{1}{2}
 (I-\Gamma).$
Thus $\Gamma$ gives rise to the decomposition $TTM= H(TM)\oplus
V(TM)$, where $H(TM):=Im \, h = Ker\, v $, $V(TM):= Im \, v=Ker \,
h$. We have
  $ J h =J, \,  h J=0,$
 $J v=0,\,  v J=J.$
The torsion $T$ of a nonlinear connection $\Gamma$ is the vector
$2$-form  on $TM$ defined by $T:=\frac{1}{2} [J,\Gamma]$.   The
curvature of a nonlinear connection $\Gamma$ is the vector $2$-form
$\Re$ on $TM$
     defined by
    $\Re:=- \frac{1}{2} [h,h].$

\begin{thm} \label{th.9a} {\em{\cite{Grifone 1}}} On a Finsler manifold $(M,L)$, there exists a unique
conservative nonlinear  connection {\em($d_{h}E=0$)} with zero
torsion. It is given by\,{\em:} \vspace{-0.3cm} $$\Gamma = [J,G]
,\vspace{-0.3cm} $$ where $G$ is the canonical spray. Such a
connection is called the canonical connection or the Barthel
connection associated to $(M,L)$.
\end{thm}

In conclusion, on a Finsler manifold $(M,L)$, there are canonically
associated  three connections\,; two of which are linear (the Cartan
connection $\nabla$ (Thm. \ref{th.1}) and the  Berwald connection
$D$ (Thm. \ref{th.1a})) and the third is nonlinear ( the Barthel
connection $\Gamma$ (Thm. \ref{th.9a})). These three connections
have the same horizontal and vertical distributions.
\bigskip
\bigskip
\bigskip
\newline
\Section{Nonlinear connections associated to a \vspace{-22pt}
regular linear connection}

 In this section we introduce two nonlinear connections associated
to a given regular connection in $\pi^{-1}(TM)$. We also introduce a
nonlinear connection naturally associated to the Cartan connection .
The relation between these three nonlinear connections and the
Barthel connection is obtained.

\begin{prop}To each  regular linear connection $D$  in  $\pi^{-1}(TM)$,
 there are associated two nonlinear connections \vspace{-0.3cm}
 $$\Gamma _{1}:= I - 2 \varphi ^{-1} o
\tilde{\varphi}, \qquad \Gamma _{2}:= 2 \beta o \rho -I
\vspace{-0.3cm},$$  where  $\tilde{\varphi}:= \gamma o K$,\,
$\varphi:=\tilde{\varphi}|_{V(TM)}$.

 Moreover, the horizontal and vertical projectors associated to $\Gamma_{1}$ and
 $\Gamma_{2}$  are given respectively
by\,{\em:}\vspace{-0.3cm}$$h_{1}=I - \varphi ^{-1} o
\tilde{\varphi},\qquad  v_{1}= \varphi ^{-1} o \tilde{\varphi}$$
 $$ \,\,\,\,h_{2} = \beta o \rho,\qquad \qquad\,\,\,\,\,v_{2}=I - \beta o \rho $$

\end{prop}
\prof Firstly, since   $ \varphi:=\tilde{\varphi}|_{V(TM)}= V(TM)
\To V(TM) $  is an isomorphism on $V(TM)$,  $\tilde{\varphi}=
\varphi$ on $V(TM)$ and $\tilde{\varphi}=0$ on $H(TM)$, it follows
that
                 $$\varphi ^{-1} o \tilde{\varphi} =\left\{
                                                            \begin{array}{ll}
                                                              id_{V(TM)} & \hbox{on V(TM)} \\
                                                              0 & \hbox{on H(TM)}
                                                            \end{array}
                                                          \right. $$
\noindent  Now,  $\Gamma _{1}$ is a vector $1$-form on $TM$
satisfying $J \Gamma _{1}=J o ( I - 2 \varphi ^{-1} o
\tilde{\varphi})=J$ and $ \Gamma _{1}J = ( I - 2 \varphi ^{-1} o
\tilde{\varphi})oJ=-J$. Hence, $\Gamma _{1}$ is a nonlinear
connection on $M$.

\par
Secondly,  $\Gamma _{2}$ is a vector $1$-form on $TM$ satisfying
   $ J \Gamma _{2} =(\gamma o \rho) o ( 2\beta o \rho-I)= J$ and
   $ \Gamma _{2} J =( 2\beta o \rho-I)o (\gamma o \rho) = - J$.
   Hence, $\Gamma _{2} $ is a nonlinear connection on $M$.\ \ $\Box$


\begin{prop} Let $(M,L)$ be a Finsler manifold. To the Cartan
connection $\nabla$, there is associated a nonlinear connection
\vspace{-0.2cm}
$$\Gamma :=\beta o \rho - \gamma o K .\vspace{-0.2cm}$$ Moreover, the horizontal and
vertical projectors of $\Gamma$ and $\nabla$ are the same.
\end{prop}

 \prof  The proof is straightforward and we
omit it.\ \ $\Box $

\begin{thm}\label{th.9} Let $(M,L)$ be a Finsler manifold. If $\,\nabla$ is the  Cartan connection
in $\pi^{-1}(TM)$, then the two nonlinear connection $\Gamma_{1}$
and $\Gamma_{2}$ associated to $\nabla$ coincide  and both equal to
\vspace{-0.2cm} $$\Gamma =\beta o \rho - \gamma o K
,\vspace{-0.2cm}$$ Moreover, $\Gamma $ coincides with the Barthel
connection associated to $(M,L)${\em:} $\Gamma=[J,G]$.
\end{thm}

\prof Since  $\nabla$ is  the Cartan connection in $\pi^{-1}(TM)$,
then $K=\gamma^{-1}$ on $V(TM)$. Therefore,  $\Gamma _{1}= I - 2
\varphi ^{-1} o \tilde{\varphi}=I-2 \gamma o K
\stackrel{(\ref{hv})}=I-2(I-\beta o \rho)= \Gamma_{2}
\stackrel{(\ref{hv})}= \beta o
\rho - \gamma o K= \Gamma .$ \\
 Moreover, the nonlinear connection $
\Gamma $
  is conservative.  In fact,  $d_{h} E(X)=i_{h}dE(X)=hX\cdot E$. But since $2E=g(\o \eta,\o
  \eta)$, by (\ref{eq.t4}), and  since $\nabla$ is metric, then
  $\frac{1}{2}d_{h}E(X)= hX \cdot g(\o \eta, \o \eta)= 2g(\o \eta, \nabla_{hX}\o \eta
  )=0.$\\
  Finally, one can show that $\Gamma$ is torsion-free. Hence, by Theorem
  \ref{th.9a}, $\Gamma $ coincides with the Barthel
   connection.\ \ $\Box$
\\
\par The above consideration enables us to express the Barthel
connection $\Gamma$ in different equivalent forms: \vspace{-0.3cm}
\begin{equation}\label{gamma}
    \Gamma= I-2
\gamma o K =2\beta o \rho - I= \beta o \rho - \gamma o K= [J,G].
\end{equation}

Note that the first three expressions of $ \Gamma$ belong to the
PB-formalism, whereas the last expression belongs to the
KG-formalism. Relations (\ref{gamma}) establishes a useful link
between the two formalisms.
\begin{lem}\label{le.1}{\em{\cite{Ali N.1}}} Under an  arbitrary change $L \To \tilde{L}$
of Finsler structures on $M$, let  the corresponding Cartan
connections $\nabla$ and $\tilde{\nabla} $ be related by
$\tilde{\nabla} _{X} \overline{Y} =\nabla _{X} \overline{Y}+\omega
(X,\overline{Y})$. If we denote \vspace{-0.3cm}
\begin{equation}\label{ABN}
  \left.
    \begin{array}{rcl}
A(\overline{X},\overline{Y})&:=&\omega(\gamma
\overline{X},\overline{Y}),\ \
B(\overline{X},\overline{Y}):=\omega(\beta
\overline{X},\overline{Y}),\ \\
N(\overline{X})&:=&B(\overline{X},\overline{\eta}),\qquad\qquad
N_{o}:= N (\overline{\eta}), \vspace{-0.3cm}
\end{array}
  \right\}
\end{equation}
\noindent then we have
\begin{description}
\item[(a)] $\omega(\overline{X},\overline{Y})=A(K{X},\overline{Y})+B(\rho {X},\overline{Y}),$

 \item[(b)] $A(\overline{X},\overline{\eta})=0,$

 \item[(c)]  $ \tilde{K}=K + N o \rho \,,  \:\:\: \tilde{\beta}=\beta -\gamma o N
  $.

\end{description}

\end{lem}
 By using Lemma \ref{le.1} and  Theorem \ref{th.9} , we have

\begin{prop}\label{pp.1} Under a change $L \To \tilde{L}$
of Finsler structures on $M$,
 the corresponding Barthel connections $\Gamma$ and $\tilde{\Gamma}
  $  are related by \vspace{-0.3cm}
  \begin{equation}
    \tilde{\Gamma} = \Gamma -2 \L,\, \text{with}\ \ \L:= \gamma o N o \rho. \vspace{-0.3cm}
\end{equation}
  Moreover, we have\ \ $ \tilde{h} =h - \L ,\ \  \tilde{v} =v + \L .$
\end{prop}

\begin{prop}\label{cor.1} The following assertion
are equivalent{\em:} \vspace{-0.2cm}
 $$(a)\, \, N=0 ,\qquad \qquad (b)\, \,     N_{o}=0,\qquad \qquad (c)\,\, \tilde{\Gamma }= \Gamma.$$
\end{prop}

 \prof   $ (a) \Longrightarrow
(b)$  is trivial. \\
 $ (b)\Longrightarrow (c)$: If $N_{o}=0$, then, by Lemma
 \ref{le.1}(b),  $\tilde{\beta}(\overline{\eta}) =\beta(\overline{\eta}) - \gamma( N_{o})=\beta(\overline{\eta}) $. Hence, $\tilde{G}=G$ and so
$\tilde{\Gamma }= \Gamma$.\\
  $ (c)\Longrightarrow (a)$: If $\tilde{\Gamma }= \Gamma$, then both $\tilde{\Gamma}$ and
   $ \Gamma$ have the same horizontal distribution. This implies that $\tilde{\beta}=\beta$. Then,
    by Lemma \ref{le.1}(b), $ \gamma o N=0$; from which $N=0$.\ \ $\Box$


\begin{rem} The map $\L$ is a vector $1$-form on $TM$ satisfying  $\L(V(TM))=0$ and $\L(T(TM))\subset
V(TM)$. Consequently,  $\L^{2}=0$ and $\L$ is thus an almost tangent
structure on $M$.
\end{rem}
\newpage
 \Section{Conformal change of Barthel connection  and\vspace{7pt}
its curvature tensor}
 \vspace{-0.3cm}
\begin{defn}
Let $(M,L)$ and $(M,\tilde {L})$ be two Finsler manifolds. The two
associated metrics $g$ and $\tilde {g}$ are said to be conformal if
there exists a positive differentiable function $\psi(x,y)$ on $TM$
such that $\tilde {g}(\overline{X},\overline{Y})=\psi (x,y)
g(\overline{X},\overline{Y})$.  Equivalently, $g$ and $\tilde{g}$
are conformal iff $\tilde {L^2}=\psi (x,y)L^2 .$ In this case, the
transformation $L\To \tilde {L}$ is said to be a conformal
transformation and the two Finsler manifold $(M,L)$ and
$(M,\tilde{L})$ are said to be conformal.
\end{defn}

\begin{prop} Let $(M,L)$ and $(M,\tilde {L})$ be two Finsler manifolds.The two
associated metrics $g$ and $\tilde {g}$ are conformal  iff the
factor of proportionality $\psi(x,y)$ is independent of the
directional argument $y$.
\end{prop}

\prof If $\psi (x,y)=\psi (x)$, then $\tilde{g}=\psi (x)g $ is
clearly a  Finsler metric on $M$, and so $g$ and  $\tilde{g}$ are
conformal.

 Conversely, let $g$ and $\tilde{g}$ be conformal. Then
 there exists a positive differentiable function $\psi (x,y)$ such
 that $\tilde{g}(\overline{X},\overline{Y})=\psi
 (x,y)g(\overline{X},\overline{Y})$.
  Setting
 $\overline{X}=\overline{Y}=\overline{\eta}$, taking into
 account the fact that $g(\overline{\eta},\overline{\eta}) =2E$, we get $\tilde {E}=\psi (x,y)E $.
Then,   $d_{J} \tilde {E}(X)=E d_{J}\psi (x,y)(X)+ \psi (x,y) d_{J}E
(X)$ for all  $ X \in \mathfrak{X}(\T M)$. But since $d_{J}E
(X)=g(\rho X,\overline{\eta})$ by (\ref{eq.t4}),
 then $\tilde{g}(\rho X,\overline{\eta})= E d_{J}\psi(x,y) (X)+
\psi (x,y)g(\rho X,\overline{\eta}).$ Therefore,  $E d_{J}\psi =0 $,
from which $ d_{J}\psi =0 $,
 and so $\psi(x,y)$ is independent of $y$.\ \ $\Box$

\vspace{7pt} From now on, we write the conformal transformation in
the form $\tilde{g}=e^{2 \sigma(x)}g $, with $\sigma(x)$ a positive
function of $x$ alone. \vspace{7pt}

\begin{defn}Let $(M,L)$ and $(M,\tilde{L})$ be two conformal Finsler
manifolds  with $\tilde{g}=e^{2\sigma(x)}g$. A geometric object $W$
is said to be  conformally invariant
    {\em(resp. conformally $ \sigma$-invariant)}
    if $\,\tilde{W}=W$ {\em(resp. $\tilde{W}=e^{2\sigma(x)}\,W$)}.
\end{defn}

\par We need the following definition for subsequent use\,:\vspace{-0.2cm}
\begin{defn}\label{grad} The vertical gradient of a function $f\in \mathfrak{F}(TM)$, denoted
$grad_{v}f$,  is the vertical vector field
 $JX$ defined by \vspace{-0.3cm}
 $$df(Y)=\bar{g}(JX,JY),\ \ \text{for all }\ Y{\!}\in
\mathfrak{X}(TM),\vspace{-0.2cm}$$ where $\bar{g}$ is the metric on
$V(TM)$ defined by~{\em\cite{Grifone 1}}\vspace{-0.3cm}
$$\bar{g}(JY,JZ)=\Omega(JY,Z),\ \ \text{for all }\ Y, Z\in \mathfrak{X}(TM).$$
\end{defn}

\par In view of Proposition  \ref{pp.1} and Proposition 2.1.5 of
\cite{Fatma}, taking the above definition  into account, we
get\vspace{-0.2cm}

\begin{thm}\label{pp.2} Let $(M,L)$ and $(M,\tilde {L})$ be  conformal Finsler manifolds with
$\tilde{g}=e^{2 \sigma(x)}g $.  The associated
  Barthel connections $\tilde{\Gamma}$ and $\Gamma$ are related by
\vspace{-0.3cm}
\begin{equation}\label{eq.1a}
  \left.
    \begin{array}{rcl}
        \tilde{\Gamma}& = &\Gamma - 2\L , \\
    \text{where} \qquad\qquad \L &:=& d\sigma \otimes \mathcal{C}+ \sigma_{1}J - d_{J}E \otimes
grad_{v}\sigma - EF= \gamma o N o \rho
 \vspace{-0.2cm}
    \end{array}
  \right\}
\end{equation}
\noindent  {\em({\it with the same notation of Proposition}
\ref{pp.1})},
 $\sigma_{1}:=d_{G}\sigma $ and  $F:=[J,grad_{v}\sigma]$.\\
      Consequently,
$\tilde{h} = h - \L $, \,\,  $\tilde{v} = v + \L $ or equivalently,
\,$\tilde{\beta} = \beta - \L o \beta  $, \,\,  $\tilde{K} = K + K o
\L $.
 \end{thm}

\prof Since $\tilde{g}=e^{2 \sigma(x)}g $, then,  using the fact
that $2E=g(\o \eta, \o \eta)$ and that $\sigma(x)$ is independent of
$y$, we get\vspace{-0.1cm}
$$\tilde{\Omega}= 2e^{2 \sigma(x)}d \sigma \wedge i_{\mathcal{C}}\Omega +e^{2 \sigma(x)}\Omega. \vspace{-0.3cm}$$
This, together with the relation  $i_{G}\Omega= -dE$, imply
that\vspace{-0.1cm}
\begin{equation}\label{spray}
   \tilde{G}= G + 2(E \ grad_{v} \sigma -\sigma_{1} \mathcal{C}).\vspace{-0.3cm}
\end{equation}
Consequently, since $\Gamma=[J,G]$, we finally get after some
manipulation\vspace{-0.1cm}
$$\tilde{\Gamma}= \Gamma -2\left \{d\sigma \otimes \mathcal{C}+ \sigma_{1}J - d_{J}E \otimes
grad_{v}\sigma - EF \right\}.\ \ \Box$$

Note that the vector form $\L$ in (\ref{eq.1a}) is expressed in the
PB-formalism by the RHS and in the KG-formalism by LHS.  Note also
that the {\it local expressions} of (\ref{eq.1a}) and (\ref{spray})
coincide with the usual local expressions found
in~\cite{Rund},\,\cite{conf. 3},...etc.

\begin{cor}In the course of the proof of Theorem \ref{pp.2}, we have
shown that{\em\,:}
\begin{description}
    \item[(a)] $\tilde{E}=e^{2 \sigma(x)}E $,

    \item[(b)]$\tilde{G}= G + 2(E \ grad_{v} \sigma -\sigma_{1}
    \mathcal{C})$,
   \item[(c)] $\tilde{\Omega}= e^{2 \sigma(x)}\Omega+2e^{2 \sigma(x)}d \sigma \wedge i_{\mathcal{C}}\Omega
    $.
\end{description}
\end{cor}
\vspace{-0.2cm}
Consequently, if the conformal change is a homothety
($\sigma$\,=\,constant), then the canonical spray $G$ is conformally
invariant and the fundamental form $\Omega$ is conformally
$\sigma-$invariant

\begin{thm}\label{pp.a2} Let $(M,L)$ and $(M,\tilde {L})$ be  conformal Finsler manifolds with
$\tilde{g}=e^{2 \sigma(x)}g $.  The curvature tensors $\tilde{\Re}$
and $\Re$ of the associated
  Barthel connections $\tilde{\Gamma}$ and $\Gamma$ are related by\vspace{-0.2cm}
\begin{equation}\label{eq.non}
\tilde{\Re}(X,Y)=\Re(X,Y)-[\L X,\L Y]-\L[hX,hY]+\mathfrak{U}_{X,Y}
\{
    v[hX,\L Y]+\L[hX,\L Y]\},\vspace{-0.2cm}
\end{equation}
 where\quad
           $\mathfrak{U}_{X,Y}\Theta(X,Y)
          =\Theta(X,Y)-\Theta(Y,X)$.
\end{thm}
\prof The proof follows from the fact that $\Re(X,Y)=-v[h X,h Y]$
(cf.\,\cite{Nabil.1}) and  that $\L[\L X,\L Y]=0$, taking Theorem
\ref{pp.2} into account.\ \ $\Box$


\Section{Conformal change of Cartan connection  and\vspace{7pt} its
curvature tensors}
 \vspace{-0.3cm}
\begin{lem}\label{le.3}Let $(M,L)$ be a Finsler manifold. Let $g$ be the Finsler
metric associated with $L$ and let $\nabla$ be the Cartan connection
 determined by the metric $g$. Then,  the following relations hold
   \begin{description}
     \item[(a)]  $ 2g(\nabla _{vX}\rho Y,\rho Z) =vX\cdot g(\rho Y,\rho Z)+
     g(\rho Y,\rho [Z,vX])+g(\rho Z,\rho [vX,Y])$.

     \item[(b)]  $ 2g(\nabla _{hX}\rho Y,\rho Z)  = hX\cdot g(\rho Y,\rho
                Z)+ hY\cdot g(\rho Z,\rho X)-hZ\cdot g(\rho X,\rho
                Y)$\\
  $  -g(\rho X,\rho [hY,hZ])+g(\rho Y,\rho
[hZ,hX])+g(\rho Z,\rho
    [hX,hY])$.
 \end{description}
\end{lem}

\prof If $\nabla$ is a metric connection in $\pi^{-1}(T M)$ with
 nonzero torsion $\textbf{T}$, one can show that  $\nabla$  is completely  determined by the relation
 \begin{equation}\label{eq}
    \begin{split}
    2g(\nabla _{X}\rho Y,\rho Z) =& X\cdot g(\rho Y,\rho
                Z)+ Y\cdot g(\rho Z,\rho X)-Z\cdot g(\rho X,\rho
                Y) \\
        & -g(\rho X,\textbf{T}(Y,Z))+g(\rho Y,\textbf{T}(Z,X))+g(\rho Z,\textbf{T}(X,Y)) \\
        & -g(\rho X,\rho [Y,Z])+g(\rho Y,\rho [Z,X])+g(\rho Z,\rho [X,Y]).
    \end{split}
 \end{equation}
 Let $\nabla$ be  the Cartan connection, then $(a)$ follows  from (\ref{eq})
by replacing $X, Y, Z$ by $vX, hY, hZ $ respectively,
  taking into account the third condition of Theorem \ref{th.1}.

   Similarly, $(b)$ follows from (\ref{eq}) by replacing
$X, Y, Z$ by $hX, hY, hZ $ respectively and using
  the second  condition of Theorem  \ref{th.1}.\ \ $\Box$
\\
\par It is worthy noting that the {\it local expressions} of (a) and (b) of the
above lemma coincide with the usual local expression found in
\,\cite{Rund},\,\cite{[Mats3]},...etc.

\begin{thm}\label{th.2} If $(M,L)$ and $(M,\tilde {L})$ are conformal Finsler
manifolds, then the associated Cartan  connections $\nabla$ and
  $\tilde{\nabla}$  are related by:\vspace{-0.2cm}
\begin{equation}\label{eq.1}
\tilde {\nabla} _{X}\overline{Y} =
    \nabla _{X}\overline{Y}+ \omega(X, \overline{Y}), \vspace{-0.3cm}
 \end{equation}
where \vspace{-0.3cm}
 \begin{equation}\label{eq.t6}
    \begin{split}
     \omega(X,\overline{Y}):=&(h X \cdot \sigma (x)) \overline{Y}
     + (\beta \overline{Y} \cdot \sigma (x)) \rho X -g(\rho X,
     \overline{Y}) \overline{P} \\
     & - T(N \overline{Y},\rho X) +  T'(\L X,\beta \overline{Y}),
     \end{split}\vspace{-0.2cm}
\end{equation}
\noindent $\overline{P}$ being a $\pi$-vector field defined by
\vspace{-0.3cm}
\begin{equation}
    g(\overline{P}, \rho Z)=hZ \cdot \sigma(x)\vspace{-0.3cm}
\end{equation}
and  $T'$ being a 2-form on $TM$, with values in $\pi^{-1}(TM)$,
defined by
 \begin{equation*}\label{eq.t5}
     g(T'(\L X,hY),\rho Z)= g({\bf T}(\L Z, h Y),\rho X).\vspace{-0.3cm}
\end{equation*}

\noindent or equivalently by\vspace{-0.4cm}

  \begin{equation}\label{eq.t5}
     g(T'(\L X,hY),\rho Z)= g(T(N \rho Z, \rho Y),\rho X).\vspace{-0.3cm}
\end{equation}
In particular,
\begin{description}

  \item[(a)]  $ \tilde {\nabla}_{\gamma \overline{X}}\overline{Y}=
    {\nabla}_{ \gamma \overline{X}}\overline{Y}$,

  \item[(b)] $\tilde {\nabla} _{\tilde{\beta}\overline{X}}\overline{Y} =
    \nabla _{\beta \overline{X}}\overline{Y}-U(\beta \overline{X}, \overline{Y}),$
\end{description}
\noindent where  \quad $U(\beta \overline{ X}, \overline{Y})=
-\omega(\beta \overline{X}, \overline{Y})+\nabla _{\L \beta
\overline{X}}\overline{Y} = - B(\overline{X}, \overline{Y})+\nabla
_{\L \beta \overline{X}}\overline{Y} $. \vspace{7pt}\\
  \end{thm}

\prof Using  Lemma \ref{le.3}(a) and Theorem \ref{pp.2}, taking into
account the fact that \\
 $\sigma=\sigma(x)$ is independent of $y$,
we get
 \begin{equation}\label{jk}
    2g(\tilde {\nabla} _{\tilde {v}X}\rho Y,\rho Z)
       = 2g(\nabla_{v X}\rho Y,\rho Z)+
      g(\rho [\L X,hY],\rho Z)+ A_{1}(X,Y,Z),
 \end{equation}
where $A_{1} $ is the 3-form on $TM$ defined by
$$A_{1}(X,Y,Z):= \L X \cdot g(\rho Y, \rho Z)+ g(\rho Y,\rho[Z,\L X]).$$
 But since $\nabla g=0$, then $A_{1}(X,Y,Z)=g(\nabla
_{\L X}\rho Y,\rho Z)+g(\textbf{T}(\L X, hY),\rho Z)$, and so
 \begin{equation}\label{eq.2}
\tilde {\nabla} _{\tilde {v}X}\rho Y= \nabla_{v X}\rho Y + \nabla
_{\L X}\rho Y
\end{equation}
Similarly, by  Lemma \ref{le.3}(b) and Theorem  \ref{pp.2}, noting
that $\rho [\L X, \L Y]=0 $,  we get
\begin{equation}\label{eq.6}
\begin{split}
   2g(\tilde {\nabla} _{\tilde {h}X}\rho Y,\rho Z) =&2g(\nabla _{hX}\rho Y,\rho Z) +
   2 hX \cdot \sigma (x) g(\rho Y, \rho Z)+ \\ &
    + 2hY \cdot \sigma (x) g(\rho X, \rho Z)
     -2hZ \cdot \sigma (x) g(\rho X, \rho Y) \\ &
      - \left \{g(\rho[hX,\L Y],\rho Z)+ g(\rho [\L X,hY],\rho Z)\right \} \\ &
      + A_{2}(X,Y,Z),
\end{split}
\end{equation}
where $ A_{2}$ is the 3-form on $TM$ defined by
\begin{equation*}\label{}
   \begin{split}
     A_{2}(X,Y,Z)   =& -\L X \cdot g(\rho Y, \rho Z) -\L Y \cdot g(\rho Z, \rho X)+\L Z \cdot g(\rho X, \rho Y)\\
       & + g(\rho X,\rho [hY,\L Z])-g(\rho Y,\rho [hZ,\L X])+ \\ & + g(\rho X,\rho [\L
       Y,hZ])-g(\rho Y,\rho [\L Z,hX]).
   \end{split}
\end{equation*}
But since $\nabla g=0$, then
\begin{equation*}\label{}
    \begin{split}
       A_{2}(X,Y,Z)  =& -g(\nabla _{\L X} \rho Y
+ \nabla _{\L Y} \rho X,\rho Z) -g(\textbf{T}(\L X,hY),\rho Z)\\
        &-g(\textbf{T}(\L Y,hX),\rho Z)+2g(\textbf{T}(\L Z,hY),\rho X).
    \end{split}
\end{equation*}
 Thus (\ref{eq.6})  reduces to
\begin{equation}\label{eq.3}
\begin{split}
  \tilde {\nabla} _{\tilde {h}X}\rho Y =&\nabla _{hX}\rho Y +
   ( hX \cdot \sigma (x)) \rho Y + (hY \cdot \sigma (x)) \rho X - g(\rho X, \rho Y) \overline{P} \\ &
      - \textbf{T}(\L X,hY) - \textbf{T}(\L Y,hX) +  T'(\L
      X,hY)-\rho[\L X,h Y].
\end{split}
\end{equation}
Now, by  (\ref{eq.2}) and (\ref{eq.3}), we get
\begin{equation}\label{eq.4}
\begin{split}
  \tilde {\nabla}_{X}\rho Y =&\nabla _{X}\rho Y +
   ( hX \cdot \sigma (x)) \rho Y + (hY \cdot \sigma (x)) \rho X - g(\rho X, \rho Y)
    \overline{P}-\rho[\L X,h Y] \\
    &- T(N \rho X,\rho Y)  - T(N \rho Y,\rho X) +  T'(\L  X,hY)+\nabla_{\L X}\rho Y.
\end{split}
\end{equation}
 Hence, the result follows from (\ref{eq.4}), making use of the
 identity   $T(N \rho X,\rho Y)=\textbf{T}(\L X,hY)= \nabla_{\L X}\rho Y-\rho[\L X,h Y] $.\ \ $\Box$
\\
\par It should be noted that the {\it local expressions} of the formulae (a) and (b) of the
above  theorem   coincide with the usual local formulae, expressing
the conformal change of  Cartan connection, found in\
\cite{Rund},\,\cite{conf. 3},\,\cite{Kitayama}...etc., where the
{\it local expression} of $\L$  plays an important  role.

\begin{rem}
For all $ X, Y \in  \mathfrak{X}(TM)$ and  $ \overline{X},
\overline{Y} \in \mathfrak{X}(\pi (M))$,
\begin{description}
\item [--]the tensor  $\omega$ satisfies the identity\,{\em:} $\omega(\gamma \overline{X},
\overline{Y})=0$.

 \item [--] the tensor $T'$
satisfies the identity\,{\em:}  $ T'(\L X,hY)= T'(\L Y,hX)$.

 \item [--]  The map $U$  satisfies the identity\,{\em:}
   $U(\beta \overline{X},\overline{\eta})=0$.

\end{description}
\end{rem}

\newpage
\par  We have already some conformal invariants and conformal
$\sigma$-invariants\,: \vspace{-0.2cm}
\begin{prop} \label{pp.6} Let the Finsler
manifolds  $(M,L)$ and $(M,\tilde{L})$ be  conformal with
$\tilde{g}=e^{2\sigma(x)}g$. Then
\begin{description}
    \item[{(a)}] If a $\pi$-tensor field $W$ of type (1,p) is conformally invariant,
    then its trace  $Tr(W)$  is conformally invariant.

    \item[(b)] The  map  $\nabla_{\gamma \overline{X}}: \mathfrak{X}(\pi (M)) \To \mathfrak{X} (\pi (M)):
  \overline{Y} \longmapsto \nabla_{\gamma \overline{X}}\overline{Y}$  is conformally
  invariant.\\
    Consequently, if $\,W$ is a conformally invariant $\pi$-tensor field, then
  so is  $\nabla_{\gamma \overline{X}}W$.

     \item[(c)] The vector $\pi$-form $\nabla\overline{X}:
     \mathfrak{X}(\pi (M)) \To \mathfrak{X} (\pi (M)):
  \overline{Y} \longmapsto \nabla_{\gamma \overline{Y}}\overline{X}$
  is  conformally invariant.

    \item[(d)] The mixed torsion $T$ of the Cartan connection is  conformally
    invariant.\\
     Consequently, $\nabla_{\gamma \overline{X}}C$ is conformally invariant.

      \item[(e)] The $\pi$-tensor $({dL\, o\, \gamma})/{L}$ is conformally
invariant{\em;} or equivalently, the tensor ${d_{J}L}/{L}$ is
conformally invariant.

    \item[(f)] The angular metric tensor $\hbar$ defined by $\hbar(\overline{X},\overline{Y})=
  g(\overline{X},\overline{Y})- \frac{1}{L^2}g(\overline{X},\overline{\eta}) g(\overline{Y},\overline{\eta})$
  is  conformally  $\sigma$-invariant.

      \item[(g)] The tensor field $\mathbb{T}$ defined by \vspace{-0.2cm}
      \begin{equation*} \mathbb{T}(\overline{X},\overline{Y}, \o Z, \o W )=
  (\nabla_{\gamma \o X}T)(\o Y, \o Z, \o W) +
  \mathfrak{S}_{\o X, \o Y, \o Z, \o W } \frac{1}{L^2}g(\o X, \o \eta) T(\o Y, \o Z, \o W)
  .
\vspace{-0.2cm}
\end{equation*}
  is   conformally $\sigma$-invariant.

\end{description}

\end{prop}

\prof  Part (a) follows from the fact that $\{\overline{E}_{i}\}$ is
an  orthonormal basis with respect to $g$ iff
$\{e^{-\sigma(x)}\overline{E}_{i}\}$ is an orthonormal basis with
respect to $\tilde{g}$. Parts  (b) and (c) follow from  Theorem
\ref{th.2}(a). Part (d) follows from the definition of $T$ and the
fact that $\rho [\gamma \o X, \tilde{\beta} \o Y]=\rho [\gamma \o X,
\beta \o Y] $. Finally, parts (e), (f) and (g) are obvious.\ \
$\Box$
\\
\par Note that the tensor $\mathbb{T}$ defined in (g) is the so-called
$T$-tensor, introduced locally by Matsumoto and
Shibata~\cite{Mats2}. Note also that some of the above invariants
globalize some Hashiguchi invariants~\cite{conf. 3}.
\par Some of the conformal invariants listed in Proposition \ref{pp.6}
can be used to characterize conformality. For example, we have

\begin{thm}\label{char}Two Finsler metrics $g$ and $\tilde{g}$ are conformal
if, and only if, \vspace{-0.3cm}
$$\frac{d_{J}\tilde{L}}{\tilde{L}}=\frac{d_{J}L}{L}.\vspace{-0.3cm}$$
\end{thm}

\prof Firstly, if $g$ and $\tilde{g}$ are conformal, then
$\frac{d_{J}\tilde{L}}{\tilde{L}}=\frac{d_{J}L}{L}$ by Proposition
\ref{pp.6}(e). Conversely, if
$\frac{d_{J}\tilde{L}}{\tilde{L}}=\frac{d_{J}L}{L}$, then
$\frac{\tilde{g}(\rho X, \o \eta)}{\tilde{L}^{2}}=\frac{g(\rho X, \o
\eta)}{L^{2}}$, by Equation(\ref{eq.t4}). Hence, \vspace{-0.4cm}
\begin{equation}\label{ee}
\tilde{g}(\rho X, \o \eta)= \phi(x,y) g(\rho X, \o \eta), \
\text{where}\ \ \phi(x,y)=\frac{\tilde{L}^2}{L^2}. \vspace{-0.3cm}
\end{equation}
For all $Y\in\mathfrak{X}(TM)$, we have $JY \cdot \tilde{g}(\rho X,
\o \eta) = JY \cdot( \phi g(\rho X, \o \eta))$; from which, since
$\tilde{\nabla}g=\nabla g=0$, \vspace{-0.3cm}
 $$ \tilde{g}(\tilde{\nabla}_{JY}\rho X,
\o \eta) + \tilde{g}(\rho X, \tilde{\nabla}_{JY}\o \eta) =
d_{J}\phi(Y) g(\rho X, \o \eta)+ \phi g(\nabla_{JY}\rho X, \o
\eta)+\phi g(\rho X,\nabla_{JY} \o \eta). \vspace{-0.3cm}$$ Using
the definition of
the Cartan torsion, Proposition \ref{pp.1} and the fact that \\
$\nabla_{JX} \o \eta = \rho X$,\ we get \vspace{-0.3cm}
$$
 \tilde{g}(\tilde{T}(\rho Y,\rho X),
\o \eta)+ \tilde{g}(\rho[J Y,h X], \o \eta) + \tilde{g}(\rho X, \rho
Y) =d_{J}\phi(x,y)(Y) g(\rho X, \o \eta) \vspace{-0.3cm}$$
$$+ \phi(x,y)g(T(\rho
Y,\rho X), \o \eta) + \phi(x,y)g(\rho[J Y,h X], \o \eta) +\phi(x,y)
g(\rho X, \rho Y).\vspace{-0.2cm}$$
 Making use of Theorem \ref{th.1}, Equation (\ref{ee}) and the
identity $T(\o X, \o \eta)=0$,  we conclude that \vspace{-0.2cm}
\begin{equation}\label{ee1}
\tilde{g}(\rho X, \rho Y)=d_{J}\phi(x,y)(Y)g(\rho X, \o \eta) +
\phi(x,y) g(\rho X, \rho Y). \vspace{-0.1cm}
\end{equation}
 Setting $X= G$ in the above equation,
we have\vspace{-0.3cm}
 $$ \tilde{g}(\rho Y, \o \eta )=d_{J}\phi(x,y)(Y)g(\o \eta, \o
\eta) + \phi(x,y) g(\rho Y, \o \eta) \vspace{-0.3cm}.$$ This,
together with (\ref{ee}), yields   $d_{J}\phi=0$. Hence, by
(\ref{ee1}), $\tilde{g}=\phi g$, where $\phi$ is \\
a (positive) function of $x$ only.\ \ $\Box$

\begin{thm}\label{th.3}Let $(M,L)$ and
$(M,\tilde{L})$ be two conformal Finsler manifolds. The curvature
tensors of the associated Cartan connections $\nabla$ and
$\widetilde{\nabla}$ are related by\,{\em:}
\begin{equation}\label{eq.5}
 \begin{split}
  \tilde{\textbf{K}}(X,Y)\overline{Z}  =& \textbf{K}(X,Y)\overline{Z}-\mathfrak{U}_{X,Y}
  \{(\nabla_{X}B)(\rho Y,\overline{Z})+ B( \rho X,B(\rho Y,\overline{Z}))\\
  &+ \frac{1}{2} B( \textbf{T}(X,Y),\overline{Z})\} \quad
  \forall\, X,Y \in \mathfrak{X}(\T M).
 \end{split}
\end{equation}
 \noindent  In particular,
 \begin{description}
  \item[(a)] $\tilde{S}(\overline{X},\overline{Y})\overline{Z}=
S(\overline{X},\overline{Y})\overline{Z}$.

  \item[(b)] $ \tilde{P}(\overline{X},\overline{Y})\overline{Z}= P(\overline{X},\overline{Y})\overline{Z}-
  V(\overline{X},\overline{Y})\overline{Z},$
  \end{description}
 \noindent where $V$ is the vector $\pi$-form  defined by\vspace{-0.2cm}
\begin{equation}
    \begin{split}
      V(\overline{X},\overline{Y})\overline{Z} = S(N\overline{X}, \overline{Y})
      \overline{Z} - (\nabla_{\gamma \overline{Y}}B)(\overline{X},\overline{Z})
     - B(T( \overline{Y},\overline{X}),\overline{Z}).
    \end{split}
\end{equation}
\begin{description}
  \item [(c)]
 $\tilde{R}(\overline{X}, \overline{Y})\overline{Z} = R(\overline{X},\overline{Y})\overline{Z}
 +H(\overline{X},\overline{Y})\overline{Z}$,
\end{description}
\noindent where $H$ is the vector $\pi$-form defined by
  \begin{equation}
    \begin{split}
    H(\overline{X}, \overline{Y})\overline{Z} =& S(N\overline{X} ,N \overline{Y})\overline{Z}
    - \mathfrak{U}_{\overline{X},\overline{Y}} \{P(\overline{X}, N \overline{Y})\overline{Z}
         +(\nabla_{\beta \overline{X}}B)( \overline{Y},\overline{Z}) \\
         & -(\nabla_{\L \beta \overline{X}}B)( \overline{Y},\overline{Z})
         +B( \overline{X},B(\overline{Y},\overline{Z}))- B(T(N
         \overline{X},\overline{Y}),\overline{Z})\}.
          \end{split}
\end{equation}

\end{thm}

\prof By Theorem \ref{th.2}, we have \vspace{-0,1cm}
\begin{align*}
   \tilde{\nabla}_{X}\tilde{\nabla}_{Y}\overline{Z}=
   \nabla_{X}\nabla_{Y}\overline{Z}+\omega( X,
   \nabla_{Y}\overline{Z})+\nabla_{X}\omega( Y,\overline{Z})+ \omega(
   X,\omega(Y,\overline{Z})).\vspace{-0,2cm}
\end{align*}
with similar expression for $
\tilde{\nabla}_{Y}\tilde{\nabla}_{X}\overline{Z}$. Moreover,
$$\tilde{\nabla}_{[X,Y]}\overline{Z}=
\nabla_{[X,Y]}\overline{Z}+\omega([X,Y],\overline{Z}).$$ The above
formulae together with the definition of the curvature tensor
$\textbf{K}$ give rise to (\ref{eq.5}). Moreover, (a) follows from
(\ref{eq.5}) by setting $X=\gamma \overline{X}$, $Y=\gamma
\overline{Y}$, noting that $\L o\gamma = 0$, $h o \L=0$ and that
 $\textbf{T}(\gamma \overline{X}, \gamma \overline{Y})=0$.
Similarly, (b) follows from the same relation by setting
$X=\tilde{\beta} \overline{X}$,  $Y=\gamma \overline{Y}$, noting
that $\L o \tilde{\beta} =\L o \beta$, $h o \L=0$ and that
$\textbf{T}(\gamma \overline{X},\L \tilde{\beta}\overline{Y})=0$.
Finally, (c) follows from the same relation by setting
$X=\tilde{\beta}\overline{X},\  Y=\tilde{\beta}\overline{Y}.\ \
\Box$
\vspace{7pt}
\par It is to be noted that the {\it local expressions} of (a), (b) and (c) of the
above  theorem   coincide with the corresponding local expressions
found in~\cite{conf. 2},\,\cite{conf. 3},\,\cite{Kitayama}...etc.
\newpage
\par In view of the above theorem, we have
\vspace{-0.15cm}
\begin{prop}\label{pp.7} The following geometric objects are
conformally invariant\,{\em:}
\begin{description}
 \item[(a)] The vertical curvature tensor $S$.

\item[(b)]
The vertical Ricci tensor  $Ric^{v}$.

\item[(c)] The scalar function $L^{2} Sc^{v}$.

\item[(d)]The $\pi$-tensor field \ $\mathbb{F}^{v}:=\{Ric^{v}-\frac{Sc^{v}\hbar}{2(n-2)}\}$.

\item[(e)]The vertical Einstein  $\pi$-tensor field
$ E^{v}:= \{Ric^{v}-\frac{Sc^{v}}{2}g\}$.

\end{description}
\end{prop}
Note that the tensor $\mathbb{F}^{v}$ in (d) is used in the
definition of the special Finsler space $S_{4}$-like.

\begin{prop}\label{pp.7} Assume that the $\pi$-tensor field  $ H$ is traceless\,{\em:} $Tr(H)=0$.
Then, the following geometric object are conformally
invariant\,{\em:}
\begin{description}

\item[(a)]
The horizontal Ricci tensor  $Ric^{h}$.

\item[(b)]  The scalar function ${L}^{2} Sc^{h}$.

\item[(c)]The $\pi$-tensor field \ $\mathbb{F}^{h}\!:=\{Ric^{h}-\frac{Sc^{h}g}{2(n-1)}\}$.

\item[(d)]The horizontal Einstein  $\pi$-tensor field
$ E^{h}:= \{Ric^{h}-\frac{Sc^{h}}{2}g\}$.

\end{description}
\end{prop}
Note that the tensor $\mathbb{F}^{h}$ in (c) is used in the
 definition of the  special Finsler space
$R_{3}$-like.

\begin{lem}\label{bracket}For all $\,\overline{X},\overline{Y}\in \mathfrak{X}(\pi(M))$, we have\,\em:

   \begin{description}
     \item[(a)] $[\gamma \overline{X},\gamma \overline{Y}]=
     \gamma(\nabla_{\gamma \overline{X}}\overline{Y}-
     \nabla_{\gamma \overline{Y}}\overline{X})$

     \item[(b)] $[\gamma \overline{X},\beta \overline{Y}]=-
     \gamma(P(\overline{Y},\overline{X})\overline{\eta}+\nabla_
     {\beta \overline{Y}}\overline{X})
     +\beta( \nabla_{\gamma \overline{X}}\overline{Y}-T(\overline{X},\overline{Y}))$

     \item[(c)] $[\beta \overline{X},\beta \overline{Y}]=
     \gamma(R(\overline{X},\overline{Y})\overline{\eta})
     + \beta(\nabla_{\beta \overline{X}}\overline{Y}-
     \nabla_{\beta \overline{Y}}\overline{X})$
 \end{description}
\vspace{-7pt} Consequently, the horizontal distribution is
completely integrable if, and only if,
$R(\overline{X},\overline{Y})\overline{\eta}=0$.

\end{lem}

\par Assume that $H(\overline{X}, \overline{Y})\overline{\eta}=0$ for all
 $\overline{X}, \overline{Y} \in \mathfrak{X}(\pi(M))$. Consequently,
$\tilde{R}(\overline{X},\overline{Y})\overline{\eta}=R(\overline{X},
\overline{Y})\overline{\eta}$, by Theorem \ref{th.3}(c). This,
together with the above lemma, give rise to the following

\begin{thm}\label{th.7} Suppose that $H(\o X, \o Y )\o \eta=0$.
The horizontal distribution with respect to $\nabla$ is completely
integrable if, and only if,  the horizontal distribution with
respect to $\tilde{\nabla}$ is completely integrable.
\end{thm}
\newpage
\par We terminate this section by the following
\vspace{-0.18cm}
\begin{thm}\label{important}Under a Finsler conformal change $\tilde{g}=e^{2 \sigma(x)}g
$, the following assertions are equivalent\,{\em:}
\begin{description}
  \item[(a)] The $\pi$-tensor field $N$ vanishes identically\,{\em:} $N=0$,
  \item[(b)]  The $\pi$-vector field $N_{o}$ vanishes identically\,{\em:} $N_{o}=0$,
  \item[(c)] The two associated Barthel connections coincide\,{\em:}
$\tilde{\Gamma}=\Gamma$,
  \item[(d)] The two associated Cartan connections coincide\,{\em:}
$\tilde{\nabla}=\nabla$,
  \item[(e)] The conformal transformation is a homothety\,{\em:}
$\sigma=$constant.
\end{description}
\end{thm}

\prof The equivalences $(a)\Longleftrightarrow
(b)\Longleftrightarrow (c)$ have been established in Proposition
\ref{cor.1}. We shall now prove the sequence of implications
$(e)\Longrightarrow (d)\Longrightarrow (a) \Longrightarrow(e)$.\\
$(e)\Longrightarrow (d)$: Let $\sigma$ be constant. Then $d
\sigma=0$, $d_{G}\sigma=0$, $grad_{v}\sigma=0$ and
$[J,grad_{v}\sigma]=0$. This implies, by (\ref{eq.1a}), that $\L =0$
(and so $N=0$). Now, putting $\sigma=$constant, $\L=0$ and $N=0$ in
(\ref{eq.t6}), we get $\omega(X, \o Y)=0$ and consequently
$\tilde{\nabla}=\nabla$.\\
$(d)\Longrightarrow(a)$: If $\tilde{\nabla}=\nabla$, then, by
(\ref{eq.1}) and (\ref{ABN}), we  have $0=\omega(\beta \o X, \o
\eta)=B(\o X, \o \eta)$; from which
$N=0$.\\
$(a)\Longrightarrow(e)$: If $N=0$, then $\L=0$ by (\ref{eq.1a}).
Now, we compute $\bar{g}(\L X, \mathcal{C})$, where $\bar{g}$ is the
metric defined in Definition \ref{grad}\,:  \\
$\bar{g}((d\sigma \otimes \mathcal{C})X, \mathcal{C})= L^{2} d \sigma(X)$.\\
$\bar{g}(\sigma_{1}J X, \mathcal{C})= d_{J}E(X) (G \cdot \sigma) $,\\
$\bar{g}((d_{J}E \otimes grad_{v}\sigma)X,\mathcal{C})= d_{J}E(X)( G \cdot \sigma)$,\\
$\bar{g}(EF(X),\mathcal{ C})= E \{
\bar{g}([JX,grad_{v}\sigma],\mathcal{C})-
\bar{g}(J[X,grad_{v}\sigma],\mathcal{C})\}=0 $,\\Substituting the
above expressions in $\bar{g}(\L X,\mathcal{C})=0$, we get
$d\sigma=0$, from which $\sigma$ is constant (provided that $M$ is
connected).\ \ $\Box$ \vspace{7pt}
\par It should be noted that some important results of
Hashiguchi~\cite{conf. 3} ({\it obtained in local coordinates}) are
thus retrieved by some parts of the above theorem.

\Section{Conformal change of Berwald connection and\vspace{7pt} its
curvature tensors}
Let   $(M,L)$ be a Finsler manifold. Let $\nabla$ and $D$ be
respectively the Cartan connection and the Berwald connection
associated to $(M,L)$. Throughout, the entities associated to the
Berwald connection will be marked by an asterisk \lq\lq\,*\,\rq\rq.

\begin{lem}\label{pp.3} {\em{\cite{Ali3}}} The horizontal maps $\beta$ and
$\beta^{*}$, associated to the Cartan connection $\nabla $
 and the Berwald connection $ D $,  coincide. Similarly,
the deflection maps $K$ and $K^{*}$ coincide.
 \end{lem}

The following result gives an explicit expression of the Berwald
connection $D$ in terms of the Cartan connection $\nabla$.

\begin{lem}\label{th.5}{\em{\cite{Ali3}}} Let $(M,L)$ be a Finsler manifold.
The Cartan  connection $\nabla $
 and the Berwald connection $ D $  are related by
  $$ D_{X}\overline{Y} = \nabla _{X}\overline{Y}
+P(\rho X, \overline{Y})\overline{\eta}- T(KX,\overline{Y})\quad
\forall\: X \in \mathfrak{X} (T M), \, \overline{Y} \in
\mathfrak{X}(\pi(M)) . \vspace{-0.2cm}$$ In particular, we have
\begin{description}
  \item[(a)] $ D_{\gamma \overline{X}}\overline{Y}=\nabla _{\gamma
  \overline{X}}\overline{Y}-T(\overline{X},\overline{Y})$.

 \item[(b)] $ D_{\beta \overline{X}}\overline{Y}=\nabla _{\beta
  \overline{X}}\overline{Y}+P(\overline{X}, \overline{Y})
  \overline{\eta}$,

\end{description}
\end{lem}

In what follows we assume that  $(M,L)$ and $(M,\tilde {L})$ are
conformal.  By Lemma \ref{pp.3} and Theorem \ref{pp.2}, we get

\begin{lem}\label{pp.4} Under a Finsler conformal change $L \longrightarrow \tilde{L}=e^{\sigma(x)}L
$, we have \vspace{-0.2cm}
 $$\tilde{h^{*}} = h^{*} - \L,\qquad  \tilde{v^{*}} = v^{*} + \L.$$
\end{lem}
\begin{thm}\label{th.4} Under a Finsler conformal change $\tilde{g}=e^{2 \sigma(x)}g $, we have
\vspace{-0.2cm}
 \begin{equation}\label{eq.2a}
    \tilde{D}_{X}\overline{Y}=D_{X}\overline{Y}+\omega^{*}(X,\overline{Y} ),\vspace{-0.2cm}
 \end{equation}

  \noindent  where\ \ $\omega^{*}(X,\overline{Y})=K([\gamma
\overline{Y}, \L]X)+D_{\L X}\overline{Y}$
\vspace{4pt}\\
 \noindent In particular, we have \vspace{-0.1cm}
\begin{description}

  \item[(a)]$\tilde{D}_{\gamma \overline{X}}\overline{Y}=D_{\gamma \overline{X}}\overline{Y}$

  \item[(b)] $\tilde{D}_{\tilde{\beta} \overline{X} }\overline{Y}=
    D_{\beta \overline{X}}\overline{Y}- \Psi(\beta
    \overline{X},\overline{Y})$,
\end{description}
 \noindent  where \,\,\,
$ \Psi(\beta \overline{X},\overline{Y})=K([\L,\gamma
\overline{Y}]\beta \o X)
  =-B^{*}( \o X, \o Y)+ D_{\L \beta \o X}\o Y.
$
\end{thm}

\prof The formula (\ref{eq.2a}) follows from Theorems \ref{th.2},
\ref{th.3} and Lemmas \ref{th.5}, \ref{pp.4}\,:
\begin{equation*}
    \begin{split}
       \tilde{D}_{X}\overline{Y}
       &=\tilde{\nabla} _{X}\overline{Y}
       +\tilde{P}(\rho X ,\overline{Y})\overline{\eta}- \tilde{T}(\tilde{K}X,
       \overline{Y}) \\
       & = \nabla _{X}\overline{Y} - U(hX,\overline{Y}) + \nabla _{\L X}\overline{Y}
       + P(\rho X ,\overline{Y})\overline{\eta}- V(\rho X ,\overline{Y})
         \overline{\eta}\\
        &   {\quad}
        -T(KX,\overline{Y}) - T(K \L X,\overline{Y})\\
& =D _{X}\overline{Y} - U(hX,\overline{Y})
       - V(\rho X ,\overline{Y})\overline{\eta}
         - T(K \L X,\overline{Y})+\nabla_{\L X}\overline{Y}\\
& =D _{X}\overline{Y}+ \nabla_{\L [hX,\gamma
\overline{Y}]}\overline{\eta} - \nabla_{[\L
hX,\gamma\overline{Y}]}\overline{\eta}+ D_{\L X}\overline{Y}\\
& =D _{X}\overline{Y}+K([\gamma \overline{Y}, \L]X)+D_{\L
X}\overline{Y}
        = D_{X}\overline{Y}+
\omega^{*}(X,\overline{Y}).
    \end{split}
\end{equation*}
 The relations  (a) and  (b) follow from (\ref{eq.2a})
by setting $X=\gamma \overline{X}$ and  $X=\tilde{\beta}
\overline{X}$\\
 respectively.\ \ $\Box$
\vspace{6pt}
\par We have to note that the {\it local expressions} of (a) and (b)
can be found in\ \cite{conf. 3},\, \cite{conf. 2},\,\cite{Rund}
...etc.\vspace{7pt}
 \par  In view of the above theorem, we have
 \vspace{-0.2cm}
\begin{prop}\label{pp.5} The $\pi$-tensor field $\omega^{*}$ has the
properties\,{\em:}
\begin{description}
      \item[(a)] $\omega^{*}(X,\rho Y)= \omega^{*}(Y,\rho X)$.

       \item[(b)]  $\omega^{*}(\gamma \overline{X},\overline{Y})= 0$.

      \item[(a)] $\omega^{*}(\tilde{\beta} \overline{X},\overline{Y})
      =\omega^{*}(\beta \overline{X},\overline{Y})$.
\end{description}
 Moreover, The map $\Psi$ has the property\,{\em:}
   $\Psi(\beta \overline{X},\overline{\eta})= 0$.
\end{prop}

\begin{thm}\label{th.6}   Under a Finsler conformal change
$\tilde{g}=e^{2 \sigma(x)}g $, we have\vspace{-0.2cm}
\begin{equation}\label{eq.7}
  \tilde{\textbf{K}^*}(X,Y)\overline{Z}  = \textbf{K}^*(X,Y)\overline{Z}-\mathfrak{U}_{X,Y}
  \{(\nabla_{X}B^{*})(\rho Y,\overline{Z})+ B^{*}(\rho X,B^{*}(\rho Y,\overline{Z})\}
\end{equation}
 for all\ $X, Y \in \mathfrak{X}(TM)$ and  $ \overline{Z} \in \mathfrak{X
  }(\pi(M)).$\\
  In particular,

\begin{description}

  \item[(a)] $\tilde{S^{*}}(\overline{X},\overline{Y})\overline{Z}=
S^{*}(\overline{X},\overline{Y})\overline{Z}=0 .$

  \item[(b)]  $ \tilde{P^{*}}(\overline{X},\overline{Y})\overline{Z}=
   P^{*}(\overline{X},\overline{Y})\overline{Z} + (D_{\gamma \overline{Y}}B^{*})(\overline{X},\overline{Z}).$

  \item[(c)]  $\tilde{R^{*}}(\overline{X}, \overline{Y})\overline{Z} =
  R^{*}(\overline{X},\overline{Y})\overline{Z}
 + H^{*}(\overline{X},\overline{Y})\overline{Z},$
 \end{description}
  \noindent  where $H^{*}$ is the vector $\pi$-form  defined by
\begin{equation*}
    \begin{split}
     H^{*}(\overline{X}, \overline{Y})\overline{Z}
       =& \mathfrak{U}_{\overline{X},\overline{Y}} \{P^{*}(\overline{Y}, N \overline{X})\overline{Z}
     - (D_{\beta \overline{X}}B^{*})(\overline{Y},\overline{Z}) \\
        &+ (D_{\L \beta \overline{X}}B^{*})(\overline{Y},\overline{Z})
     -B^{*}( \overline{X},B^{*}(\overline{Y},\overline{Z}))\}.
    \end{split}
\end{equation*}

\end{thm}

\prof  Equation (\ref{eq.7}) follows from the definition of the
curvature tensor $\textbf{K}^*$,
 together with Theorem \ref{th.4} and Lemma \ref{pp.4}.
  Part (a) follows by setting $X= \gamma \overline{X}$ and $Y=\gamma \overline{Y}$
   in (\ref{eq.7}), taking  Proposition \ref{pp.5} into account.
 Relation (b) follows by setting $X= \tilde{\beta }\overline{X}$ and $Y=\gamma \overline{Y}$
   in (\ref{eq.7}), taking   Lemma \ref{pp.4} and Proposition \ref{pp.5} into
   account.
   Relation (c) follows by setting $X= \tilde{\beta} \overline{X}$ and $Y=\tilde{\beta} \overline{Y}$
   in  the same equation.\ \ $\Box$
\\
\par The {\it local expressions} of (b) and (c) of the
above  theorem   are the same as those found in\ \cite{conf.
2},\,\cite{Rund},\,\cite{conf. 3},...etc.


\begin{prop} Let the Finsler
manifolds  $(M,L)$ and $(M,\tilde{L})$ be  conformal.
\par The following geometric objects are conformally invariant\,{\em:}\vspace{-0.3cm}
\begin{description}
\item[(a)] The  map  $D_{\gamma \overline{X}}:
\mathfrak{X}(\pi (M)) \To \mathfrak{X} (\pi (M)):
  \overline{Y} \longmapsto D_{\gamma \overline{X}}\overline{Y}$.
  \item[(b)]The vector $\pi$-form $D\overline{X}:
     \mathfrak{X}(\pi (M)) \To \mathfrak{X} (\pi (M)):
  \overline{Y} \longmapsto D_{\gamma \overline{Y}}\overline{X}$.
 \end{description}
\par   Given that the $\pi$-form $B^{*}$ is vertically
parallel\,{\em:} $D_{\gamma\overline{X}}B^{*}=0$,
then\vspace{-0.3cm}
\begin{description}
\item[(c)]The mixed curvature tensor $P^*$ is conformally
invariant.
\end{description}
\par  Given that the $\pi$-tensor field $H^{*}$ is traceless{\,\em:} $
Tr(H^{*})=0$, then the following geometric objects are conformally
invariant\,{\em:}\vspace{-0.3cm}
\begin{description}
\item[(d)]The horizontal Ricci tensor  ${Ric^{*}}^{h}$.
\item[(e)]The scalar function $L^{2}{Sc^{*}}^{h}$.
\item[(f)]The horizontal Einstein  $\pi$-tensor field $
{E^{*}}^{h}:= \{{Ric^{*}}^{h}-\frac{{Sc^{*}}^{h}}{2}g\}$.
\end{description}
\end{prop}

\Section{Application: Geodesics and Jacobi fields} In this section
we present an application of some of the obtained results.
\par Let $c: I \longrightarrow M$ be a regular curve in $M$. The
canonical lift of $c$ is the curve $\hat{c}$ in $\T M$ defined by $
\hat{c}:t\longmapsto dc\diagup dt$. The lift of a vector field $X
\in \mathfrak{X}(M)$ along $c$ is the $\pi$-vector field along
$\hat{c}$ defined by $\overline{X}:\hat{c}(t)\longmapsto
(\hat{c}(t),X(c(t)))$. In
 particular, the velocity vector field $dc \diagup dt$ along $c$ is
 lifted to the $\pi$-vector field $\overline{d}c\diagup dt:=(dc\diagup dt,dc\diagup
 dt)$ along $\hat{c}$.
 Clearly, $\rho (d\hat{c}\diagup dt)=\overline{d}c\diagup dt=\overline{\eta}\mid_{\hat{c}(t)}$.
  A vector field $X$ along a regular curve $c$ in $M$ is parallel
along $c$ with respect to the connection $\nabla$ (or is
$\nabla$-parallel along $c$) if $D\overline{X}\diagup
dt=\nabla_{d\hat{c} \diagup dt}\overline{X}=0$, where $D\diagup dt$
\,is the covariant derivative operator associated  with $\nabla$,
along $\hat {c}$. A regular curve $c$ in $M$ is a geodesic if the
$\pi$-vector field $\frac{D}{dt}(\frac{\overline{d}c}{dt})$ vanishes
identically. In this case, the vector field $d\hat{c}\diagup dt$
along $\hat{c}$ is horizontal with respect to $\nabla$.

 In what follows we assume that  $(M,L)$ and
$(M,\tilde {L})$ are conformal Finsler manifolds with
$\tilde{g}=e^{2 \sigma(x)}g $ . From Theorem \ref{th.2}, the
associated Cartan connections $\nabla$ and $\tilde{\nabla}$ are
related by
 \vspace{-0.1cm}
\begin{equation}\label{eq.9}
   \tilde {\nabla}_{X}\overline{Y} =\nabla _{X}\overline{Y} + \omega
(X,\overline{Y} ),\vspace{-0.2cm}
\end{equation}
where $\omega$ is given by (\ref{eq.t6}).

 Let $D/dt$ and $\tilde{D}/dt$ be the associated
covariant operators (along a curve $\hat{c}$ in $\T M$)  with
respect to the Cartan connections $\nabla$ and $\tilde{\nabla}$
respectively. Then for every $\pi$-vector field $\overline{X}$ along
$\hat{c}$, taking (\ref{eq.9}) into account, we get\vspace{-0.2cm}
\begin{equation}\label{eq.t1}
   \tilde{D} \overline{X}/dt= D \overline{X}/dt + \omega
(d \hat{c}/dt, \overline{X}) .\vspace{-0.2cm}
\end{equation}
 A direct consequence of (\ref{eq.t1}) is the following

\begin{lem} Let $c$ be a regular  curve in $ M$ and $X$ a vector field along $c$.
If $X$ is  $\nabla$-parallel {\em({\it resp.
$\tilde{\nabla}$-parallel})} along $c$, then a necessary and
sufficient condition for  $X$ to be $\tilde{\nabla}$-parallel
{\em({\it resp. $\nabla$-parallel})} along $c$ is that $\omega (d
\hat{c}/dt, \overline{X})=0$.
\end{lem}

\begin{thm} A necessary and sufficient condition for a geodesic $c$ in $(M,L)$ {\em(resp. $(M,
\tilde{L})$)} to be a geodesic in $(M,\tilde{L})$ {\em(resp. $(M,L
)$)} is that $B(\overline{\vartheta},\overline{\vartheta})=0$, where
$\overline{\vartheta}:=\overline{\eta}|_{\hat{c}(t)}$ and $B$ is the
$\pi$-tensor field defined by {\em(\ref{ABN})}.
\end{thm}
 \prof Let $c$ be a regular curve in $M$ and let $\hat{c}$ be its canonical lift
 to $\T M$. Then, using Equation (\ref{eq.t1}) and the fact that
 $A(\overline{X}, \overline{\eta})=0$ (Lemma \ref{le.1}(a)), we
obtain \vspace{-0.2cm}
\begin{equation*}
  \tilde{D} \overline{\vartheta}/dt  = D \overline{\vartheta}/dt + \omega
( \beta \rho  (d \hat{c}/dt), \overline{\vartheta})
         = D \overline{\vartheta}/dt + \omega
(\beta \overline{\vartheta}, \overline{\vartheta}) = D
\overline{\vartheta}/dt + B(\overline{\vartheta},
\overline{\vartheta}).\vspace{-0.2cm}
\end{equation*}
The result follows from the last relation noting that $c$ is
geodesic in $(M,L)$ (resp. $(M, \tilde{L})$) iff
$D\overline{\vartheta}/dt=0$ (resp.
$\tilde{D}\overline{\vartheta}/dt=0$).\ \ $\Box$
\begin{defn}\label{def.3} {\em{\cite{Ali N.1}}} A vector field $\xi{\!\!}\in{\!\!}\mathfrak{X}(M)$
along a geodesic $c$ in $M$ is called a Jacobi field with respect to
a regular connection  $\nabla$ in $\pi^{-1}(TM)$ if it satisfies the
Jacobi differential equation \vspace{-0.4cm}
$$D^{2}\overline{\xi}/dt^{2}+
R(\overline{\vartheta},\overline{\xi})\overline{\vartheta}=0,\vspace{-0.2cm}$$
where
  $R$ is the $h$-curvature of $\,\nabla$,
   $\,\overline{\xi}$ is the lift
of $\,\xi$ along $c$ and
$\,\overline{\vartheta}=\overline{\eta}|_{\hat{c}}$.
\end{defn}

\begin{thm}\label{th.8}Let $c$ be a geodesic in $M$ and $\overline{\vartheta}=\o \eta|_{\hat{c}(t)}$. Assume that
 $H(\overline{\vartheta}, \overline{X})\overline{\vartheta}=0$ and that the $\pi$-tensor field
 $i_{\o \vartheta} B$ vanishes.   A vector field $\xi\in\mathfrak{X}(M)$
along $c$ is a Jacobi field with respect to $\nabla$ if, and only
if, it is a Jacobi field with respect to $\tilde{\nabla}$.
\end{thm}

\prof It should firstly be noted that a geodesic $c$ in $(M,L)$ is
also a geodesic in $(M,\tilde{L})$ since $i_{\o \vartheta} B=0$. By
hypothesis, we have\vspace{-0.2cm}
\begin{equation}\label{eq.t2}
   \tilde{R}(\overline{\vartheta}, \overline{X})\overline{\vartheta}=
R(\overline{\vartheta}, \overline{X})\overline{\vartheta} \quad
 \forall \,\, \overline{X}\in \mathfrak{X}(\pi(M)).\vspace{-0.2cm}
\end{equation}

Now, let $\overline{\xi}$ be the lift of the  vector field $\xi$
along $c$. Putting $\overline{X}=\overline{\xi}$ in (\ref{eq.t1})
and noting that $\beta o \rho  +\gamma o K = id_{\mathfrak{X}(\T
M)}$, we get\vspace{-0.2cm}
\begin{equation*}
    \begin{split}
      \tilde{D} \overline{\xi}/dt & =D \overline{\xi}/dt + \omega
(\beta \rho (d \hat{c}/dt), \overline{\xi}) + \omega(\gamma K(d \hat{c}/dt), \overline{\xi})  \\
       & = D \overline{\xi}/dt + B(\overline{\vartheta}, \overline{\xi}) + A(K(d \hat{c}/dt), \overline{\xi}).
    \end{split}
    \end{equation*}
Moreover,  since $c$ is a geodesic and
 since   $ B(\o \vartheta, \o X)=0$,  it follows that
$\tilde{D} \overline{\xi}/dt  =D \overline{\xi}/dt$.
Consequently,\vspace{-0.2cm}
\begin{equation}\label{eq.t3}
\tilde{D}^{2} \overline{\xi}/dt^{2}  =D^{2}
\overline{\xi}/dt^{2}.\vspace{-0.1cm}
\end{equation}
 According to Definition \ref{def.3},  the result follows from (\ref{eq.t2}) and
 (\ref{eq.t3}).\ \ $\Box$

\vspace{20pt}
 {\noindent\Large\textbf{Concluding remarks}}\\

\vspace{-0.3cm}
   \noindent {$\bullet$}\  A global theory of conformal Finsler geometry is
    established. Some known results are generalized and several new
    results are obtained.\\
    \textbf{$\bullet$}\ It is shown that the Pull-back formalism and the
    Klein-Grifone formalism are not alternatives but rather
    complementary.\\
    {$\bullet$}\  Although our treatment is entirely global, the local
    expressions of the obtained results, when calculated, coincide
    with the known classical local results (See the Appendix).\\
    {$\bullet$}\ The conformal change of different types of special Finsler
     spaces is not treated in the present work. It merits a separate
     study that we are currently in the process of preparing, and it
     will be the object of a forthcoming paper.\\
   {$\bullet$}\ The most important and well known connections in Finsler
    geometry are the Cartan, Berwald and Barthel connections. However,
     there are other connections  of  particular importance, such as
     Hashiguchi and Chern (Rund) connections. Such connections are not treated
     here.
     We are currently investigating these connections (and their conformal transforms)
     from a global standpoint, in the hope to build, with the material we have,
     a global theory of Finsler
     geometry as complete as possible.
\vspace{24pt}
\begin{center}
{\bf\Large{Appendix. Local formulae}}\end{center}
 \vspace{3pt}
 \par For the sake of completeness, we present in this appendix
 a brief and concise survey of the local expressions of the most
 important geometric objects treated in the paper.
  \par Let $(U,(x^{i}))$ be  a system  of local coordinates on
 $M$ and $(\pi^{-1}(U),(x^i,y^i))$ the associated system of local coordinates on $TM$.
 We use the following notations\,:\\
 $(\pa_{i}):=(\frac{\pa}{\pa x^i})$: the natural basis of $T_{x}M,\, x\in
 M$,\\
 $(\paa_{i}):=(\frac{\pa}{\pa y^i})$: the natural basis of $V_{u}(\T M),\, u\in
 \T M$,\\
$(\pa_{i},\paa_{i})$: the natural basis of $T_{u}(\T M)$,\\
$(\o \pa_{i} )$: the natural basis of the fiber over $u$ in $\p$
($\o \pa_{i} $ is the lift of $\pa_{i}$ at $u$).
\par To  a Finsler manifold $(M,L)$, we associate the geometric
objects\,:\\
$g_{ij}:= \frac{1}{2} \paa_{i} \paa_{j}L^2= \paa_{i} \paa_{j}E$: the
Finsler metric tensor,\\
$G^h$: the components of the canonical spray,\\
$G^{h}_{i}:=\paa_{i}G^h$,\\
$G^{h}_{ij}:=\paa_{j}G^h_{i}=\paa_{j}\paa_{i}G^h$,\\
$(\delta_{i}):=(\pa_{i}-G^{h}_{i}\paa_{h})$: the basis of $H_{u}(\T
M)$ adapted to $G^{h}_{i}$,\\
$(\delta_{i}, \paa_{i})$: the basis of $T_{u}(\T M)=H_{u}(\T
M)\oplus V_{u}(\T M)$ adapted to $G^{h}_{i}$.
\par We have\,:\\
$\gamma(\o \pa_{i})=\paa_{i}$,\\
$\rho(\pa_{i})=\o \pa_{i}$, \ \ $\rho(\paa_{i})=0$, \ \
$\rho(\delta_{i})=\o \pa_{i}$,\\
$\beta(\o \pa_{i})=\delta_{i}$,\\
$J(\pa_{i})= \paa_{i}$, \ \ $J(\paa_{i})=0$, \ \
$J(\delta_{i})= \paa_{i}$,\\
$h= \beta o \rho= dx^{i} \otimes \pa_{i}- G^{i}_{j}\, dx^{j} \otimes
\paa_{i}$,  \ \    $v=\gamma o K=dy^{i} \otimes \paa_{i}+
G^{i}_{j}\, dx^{j} \otimes \paa_{i} $.
\par We define\,:\\
\noindent $\gamma^{h}_{ij}:= \frac{1}{2}\,g^{h\ell}(\pa_{i}\,g_{\ell
j}+\pa_{j}\,g_{i\ell }- \pa_{\ell}\,g_{i j}
 )$,\\
  $C^{h}_{ij}:= \frac{1}{2}\,g^{h\ell}(\paa_{i}\,g_{\ell j}+\paa_{j}\,g_{i\ell }-
   \paa_{\ell}\,g_{i j})=
  \frac{1}{2}\,g^{h\ell}\,\paa_{i}\,g_{\ell j}$ \  (cf. Lemma \ref{le.3}(a)
  ),\\
  $\Gamma^{h}_{ij}:= \frac{1}{2}\,g^{h\ell}(\delta_{i}\,g_{\ell j}+\delta_{j}\,g_{i\ell }-
  \delta_{\ell}\,g_{i j})$ \  (cf. Lemma \ref{le.3}(b)).\\
   Then, we have\,:\\
   $\bullet$ The canonical spray $G$: $G^h=\frac{1}{2}\,\gamma^{h}_{ij}\,y^i
   y^j$.\\
  $\bullet$ The Barthel connection
    $\Gamma$: $G^{h}_{i}= \paa_{i}G^{h}=G^{h}_{ij}\,y^i
  =\Gamma^{h}_{ij}\,y^j=\Gamma^{h}_{io}$ \  (cf.  Equation ( \ref{gamma})).\\
 $\bullet$ The Cartan connection  $C\Gamma$:
  $(G^{h}_{i},\ \Gamma^{h}_{ij},\ C^{h}_{ij})$.\\
$\bullet$ The Berwald connection  $B\Gamma$: $(G^{h}_{i},\ G^{h}_{ij},\ 0)$.\\
We also have $G^{h}_{ij}= \Gamma^{h}_{ij}+ C^{h}_{ij}{_{|k}}\,y^k=
\Gamma^{h}_{ij}+ C^{h}_{ij}{_{|o}}$, where the the stroke \lq \lq\,
$|$ \,\rq \rq \  denotes the  horizontal Cartan covariant derivative
(cf. Lemma \ref{th.5}(b)).\\

\par Under a conformal change $g_{ij} \To \tilde {g}_{ij}= e^{2 \sigma(x)}g_{ij}$, we have
the following expressions for the relationships between various
geometric objects and their conformal transforms\,:

\vspace{7pt}
 \noindent \textit{\textbf{$\bullet$\ \ Canonical spray}}

\vspace{7pt}
 \noindent $ \tilde{G}^{h}= G^{h}- B^{h}$,\\
  where\ \  $B^{h}:= (E g^{hj}-y^hy^j)\sigma_{j}$; \, $\sigma_{j}:= \pa_{j} \sigma$.  \\
 (local expression of Equation ( \ref{spray}) ).

\vspace{7pt}
 \noindent \textit{\textbf{$\bullet$\ \ Barthel connection}}

\vspace{7pt}
  \noindent $ \tilde{G}^{h}_{j}= G^{h}_{j}- B^{h}_{j}$,\\
  where \  $B^{h}_{j}:=\paa_{j}B^h =  y_{j} \sigma^{h}-\delta^{h}_{j}\sigma_{o}-y^{h}\sigma_{j}-L^2
  C^{h}_{j}$,\,  $ \sigma_{o}:=\sigma_{i}y^{i}$,\, $\sigma^{h}:=g^{hj}\sigma_{j}$,\, $C^{i}_{j}:=C^{ir}_{j}
  \sigma_{r}$ and  $C^{ir}_{j}:=g^{rk}C^{i}_{jk}$.   \\
  (local expression of Equation ( \ref{eq.1a})).

\vspace{7pt}

   \noindent \textit{\textbf{$\bullet$\ \ Barthel  curvature tensor}}

\vspace{7pt}
 \noindent $
 \tilde{\Re}^{h}_{ij}=\Re^{h}_{ij}+H^{h}_{ij}$,\\
   where\ \  $H^{h}_{ij}:=-\mathfrak{U}_{ij}\{B^h_{i|j}+(B^h_{im}-P^h_{im})  B^m_{j}
     \};\ \ P^h_{im}:=C^h_{im|o} $.\\
 (local expression of  Equation (\ref{eq.non}) ).

\vspace{7pt} \noindent \textit{\textbf{$\bullet$\ \ Cartan
connection}}

\vspace{7pt}
  \noindent  $\textbf{\bf--}\ \  \tilde{C }^{h}_{ij}= C^{h}_{ij}$.\\
  $\textbf{\bf--}\ \ \tilde{\Gamma}^{h}_{ij}= \Gamma^{h}_{ij}- U^{h}_{ij},$\\
  where\ \ $U^{h}_{ij}:= g_{ij}\sigma^{h}- \delta^{h}_{i} \sigma_{j}-\delta^{h}_{j}
  \sigma_{i}- C^{h}_{im}B^{m}_{j}-C^{h}_{jm}B^{m}_{i}+g^{hr} C_{ijm}
  B^{m}_{r}$.\\
    (local expressions of  Theorem \ref{th.2}(a), (b) ).
\pagebreak

\vspace{7pt}
 \noindent \textit{\textbf{$\bullet$\ \ Cartan Curvature tensors}}

\vspace{7pt}
\noindent $\textbf{\bf--}\ \ \tilde{S}^{h}_{kij}= S^{h}_{kij} $.\\
$\textbf{\bf--}\ \ \tilde{P}^{h}_{kij}= P^{h}_{kij} - V^{h}_{kij} $, \\
  where \ \ $ V^{h}_{kij}:= 2B ^{m}_{i}S^{h}_{kjm}+\dot{\partial}_{j}A^{h}_{ki}-
 U^{h}_{im}C^{m}_{kj}+U^{m}_{ki}C^{h}_{jm}$;  $A^{h}_{ij}= U^{h}_{ij}+
 C^{h}_{im}B^{m}_{j}.$\\
  $\textbf{\bf--}\ \ \tilde{R}^{h}_{kij}= R^{h}_{kij} +H ^{h}_{kij} $,\\
    where\ \ $H^{h}_{kij}:= 2 S^{h}_{kml}B^{m}_{i}B^{l}_{j}- \mathfrak{U}_{ij}
\{A^{h}_{ki | j }+ B^m_{j}\dot{\partial}_{m}A^{h}_{ki}+
U^{m}_{kj}U^{h}_{im}-B^{m}_{j}P^{h}_{kim}\}. \ \ $\\
 (The local expressions of Theorem \ref{th.3}(a), (b), (c)).

\vspace{7pt}
 \noindent \textit{\textbf{$\bullet$\ \ Berwald connection}}

\vspace{7pt}
 \noindent $\textbf{\bf--}\ \ \tilde{C}^{*h}_{\ \,ij}=C^{*h}_{\ \,ij}= 0 $, \\
        where\ \ $D_{\dot{\partial_{i}}}\overline{\partial}_{j}=: C^{*h}_{\ \,ij}
       \,\, \overline{\partial}_{h}$.\\
         $\textbf{\bf--}\ \ \tilde{G}^{h}_{ij}=G^{h}_{ij} - \Psi^{h}_{ij} $,\\
          where\ \  $D_{e_{i}}\overline{\partial}_{j}=:
          G^{h}_{ij} \,\overline{\partial}_{h}$, \,\,
     $\Psi(e_{i},\overline{\partial}_{j})=: \Psi^{h}_{ij}\, \overline{\partial}_{h}
     =\dot{\partial_{j}} B^{h}_{i} \,\overline{\partial}_{h}=:B^{h}_{ij}
    \, \overline{\partial}_{h}$.\\
     (The local expressions of Theorem \ref{th.4}(a), (b) ).


\vspace{7pt}
 \noindent \textit{\textbf{$\bullet$\ \ Berwald Curvature tensors}}

\vspace{7pt}
\noindent  $\textbf{\bf--}\ \ \tilde{S}^{*h}_{\ \,kij}= S^{*h}_{\ \,kij}=0 $.\\
  $\textbf{\bf--}\ \ \tilde{P}^{*h}_{\ \,kij}= P^{*h}_{\ \,kij} - V^{*h}_{\ \,kij} $,\\
   where \  \ $ V^{*h}_{\ \,kij}:= \dot{\partial_{j}}B^{h}_{ki} $\ \
   (note that  $ P^{*h}_{\ \,kij}= \dot{\partial_{j}}G^{h}_{ki}$).\\
      $\textbf{\bf--}\ \ \tilde{R}^{*h}_{\ \,kij}= R^{*h}_{\ \,kij} +H ^{*h}_{\ \,kij} $,\\
     where\ \ $H^{*h}_{\ \,kij}:= \mathfrak{U}_{ij} \{
(\dot{\partial}_{m}G^{h}_{ik})B^{m}_{j}-B^{h}_{ik(j)}-(\dot{\partial_{m}}
B^{h}_{ik})B^{m}_{j}- B^{h}_{im}B^{m}_{kj}\} $.\\
 The parentheses \lq\lq\,$(\,)$\,\rq\rq  \ denote the horizontal Berwald  covariant
derivative.\\
  (The local expression of Theorem \ref{th.6}(a), (b), (c)).


\end{document}